\documentclass[11pt,oneside]{amsart}

\usepackage{enumitem,kantlipsum}
\usepackage[table]{xcolor} 
\usepackage{booktabs} 
\usepackage{adjustbox} 
\usepackage{cancel}
\usepackage{epigraph}
\usepackage[a4paper]{geometry}
\usepackage[dvipsnames]{xcolor}
\usepackage{graphicx}
\usepackage{amssymb}
\usepackage{amssymb}
\RequirePackage{doi}
\usepackage{hyperref}
\usepackage{graphicx}
\usepackage{forest}
\usepackage{tikz}
\usetikzlibrary{arrows.meta}
\usepackage{tikz-cd}
\usepackage{stmaryrd}
\usepackage{enumitem}
\usepackage{comment} 
\usepackage{newtxmath}
\usepackage{nicefrac} 
\usepackage{xstring}
\usepackage{subcaption} 
\usepackage{multirow}
\usepackage{tabularx}

\usepackage{xcolor}

\def\SS{\mathbb{S}}

\DeclareMathOperator{\CT}{PriorityTree}
\DeclareMathOperator{\PT}{ParkingTree}
\DeclareMathOperator{\pv}{pv}
\DeclareMathOperator{\pt}{pt}
\DeclareMathOperator{\diff}{diff}

\DeclareMathOperator{\ones}{ones}
\DeclareMathOperator{\lucky}{lucky}

\DeclareMathOperator{\leaves}{leaves}
\DeclareMathOperator{\dis}{dis}

\DeclareMathOperator{\rec}{rec}

\DeclareMathOperator{\Rec}{Rec}

\DeclareMathOperator{\chseq}{chseq}
\DeclareMathOperator{\mult}{mult}
\DeclareMathOperator{\psa}{psa}

\DeclareMathOperator{\probes}{probes}
\DeclareMathOperator{\wait}{wait}
\DeclareMathOperator{\len}{len}
\DeclareMathOperator{\ord}{ord}
\DeclareMathOperator{\lpr}{lpr}
\DeclareMathOperator{\inv}{inv}

\newcommand{\Pref}{\pi_T}

\newcommand{\pfunction}[2][matrix]{%
  \begingroup
    \setlength{\arraycolsep}{3pt}
    \StrSubstitute{#2}{;}{\\}[\pfunc@tmp]%
    \StrSubstitute{\pfunc@tmp}{,}{ & }[\pfunc@tmp]%
    \begin{#1}
      \pfunc@tmp
    \end{#1}
  \endgroup
}

\definecolor{USred}{cmyk}{0,1.00,0.65,0.34}
\renewcommand{\emph}[1]{{\textcolor{USred}{\em #1}}}

\newcommand{\memph}[1]{{\textcolor{USred}{#1}}}

\hypersetup{
  colorlinks,
  citecolor=USred,
  linkcolor=USred,
  urlcolor=USred,
  allcolors=USred}

\usepackage{amsthm}
\theoremstyle{definition}
\newtheorem{de}{Definition}[section]
\theoremstyle{plain}
\newtheorem{thm}[de]{Theorem}

\newtheorem{lem}[de]{Lemma}

\newtheorem{cor}[de]{Corollary}
\newtheorem{ex}[de]{Example}
\theoremstyle{remark}
\newtheorem{remark}[de]{Remark}

\theoremstyle{definition}

\title{
On weary drivers, records of trees, and parking functions.}
\author{Adrián Lillo}
\author{Mercedes Rosas}
\author{Stefan Trandafir}

\address[A. Lillo, M. Rosas]{Departamento de Álgebra, Universidad de Sevilla}
\email[A. Lillo]{alillo@us.es}
\email[M. Rosas]{mrosas@us.es}

\address[S. Trandafir]{Departamento de Física Aplicada II, Universidad de Sevilla}
\email[S. Trandafir]{strandafir@us.es}

\begin{document}

\begin{abstract}
This work builds on the notion of record of rooted trees. We provide an alternative definition of parking functions,  derive from it a record-preserving bijection between rooted trees and parking functions, and establish a join equidistribution result between a 5-tuple of statistics on rooted trees and a corresponding 5-tuple of statistics on parking functions. Some enumerative questions are also considered.
\end{abstract}

\maketitle

\tikzset{
    record/.style={
        circle,
        draw=black,
        fill=black,
        inner sep=1.5pt,
        label={#1},
        solid
    },
    record/.default={},
    non-record/.style={
        circle,
        draw = black,
        fill = white,
        inner sep=1.5pt,
        solid
    },
    every label/.style={
            font=\tiny
        }
}

\definecolor{CB.red}{HTML}{D83034}
\definecolor{CB.darkblue}{HTML}{003A7D}
\definecolor{CB.blue}{HTML}{008DFF}
\definecolor{CB.pink}{HTML}{003A7D}
\definecolor{CB.orange}{HTML}{FF9D3A}
\definecolor{CB.green}{rgb}{0.0,0.4,0.0}
 

\vspace{-.4cm}
\epigraph{\em
        All these cars trying to overtake each other, only to end up having to find a place to park. 
}

\section{Introduction} \label{sec:intro}

\noindent
{\bf Act I: }
Every morning, cars travel on a long one-way road, each driver with a different favorite speed.
The road is long enough that any faster car will inevitably catch up to a slower one. Overtaking is impossible. Days when the slowest car arrives first on the road are minor calamities. A traffic jam is no way to start the day.   Most mornings aren’t quite that grim, but still, it feels like every car is doomed to crawl behind someone slower. \vspace{.3cm}

 What does an average morning look like?  Assume that cars arrive from the right and are ranked from 1 as the fastest to $n$ as the slowest. The order in which the cars enter the road defines a permutation. For emphasis, we sometimes insert dots to indicate the separation between different clusters.  
 
 A cluster forms when a car gets stuck behind a cluster of slower cars. The first car to arrive always starts a cluster. When the second car arrives, it is either the slowest (2) and starts a new cluster, resulting in permutation $1 \cdot 2$,  or it is the fastest (1) and does not, as in $21$. The second car starts a new cluster half the time.
 Consider three cars on the road. If the slowest car ($3$) arrives last, it cannot catch up and starts a new cluster, resulting in $2  1 \cdot3$ or $1\cdot 2 \cdot3$. If the medium-speed car ($2$) arrives last, it merges into an existing cluster, resulting in  $3 12$ or $1 \cdot 3 2$. If the fastest car ($1$) arrives last, it becomes stuck behind the others, resulting in $2\cdot 31$ or $321$. In one-third of the cases, the third car forms a new cluster.
 A similar analysis holds for any number of cars. Only when the slowest car arrives last, a new cluster is formed. 
 Therefore, when $n$ cars drive down the road the average number of clusters is
$
1 + \frac{1}{2} + \frac{1}{3} + \ldots  + \frac{1}{n}$.

Probabilists have explored this situation from multiple angles \cite{Renyi, Glick, Wilf, Kortchemski}. Combinatorialists, on the other hand, will recognize  Foata's fundamental transformation \cite{Foata} behind our analysis. Cars at the head of the clusters are the records (left-to-right maxima) of the permutation recording the order of arrival of the cars.

\vspace{.3cm}

\noindent{\bf Act II: }
The situation has changed. Now there are junctions where smaller roads merge into larger ones. Yet, overtaking remains impossible. A new scenario, but cars still inevitably find themselves trapped behind slower drivers. 

\vspace{.3cm}

Instead of using permutations and their records to analyze this problem, we use finite labelled rooted trees. 
Cars are represented by the nodes of the tree, with edges corresponding to the roads. Node labels indicate the speed ranges of the cars. As before, we assume distances are large enough for any faster car to catch up with a slower one right ahead.

A node of a rooted tree is said to be a \emph{record} if it has the largest label along the unique path from itself to the root.  Building on this notion, we introduce 
in Section~\ref{suse:record} the \emph{record decomposition} of a rooted tree, closely related to an encoding of trees independently developed by Kreweras and Moszkowski \cite{Kreweras-Moszkowski}, and Picciotto \cite{Picciotto}.
Tree records appear implicitly in the work of A. Meir and J.~W. Moon \cite{MM} on a random process where, given a rooted tree \( T \), an edge of $T$ is randomly selected and deleted. Among the two resulting trees, the one containing the root is retained, while the other is discarded. This process continues until only the root remains. The connection between this process and the notion of records was first described in \cite{Janson}. 
 
 Although labelled trees are among the most extensively studied structures in combinatorics, combinatorialists seem to have paid little attention to the study of records of rooted trees. There are exceptions; Gessel and Seo \cite{Gessel-Seo} have used the dual version of this notion in their hook length formulas for trees.
  In contrast, related topics such as tree descents and tree inversions have attracted considerable interest \cite{M-R, G-S-Y, Yan}. 
      In our forthcoming work \cite{LRT-GF-records}, we introduce the record numbers, which enumerate rooted trees by their number of records.
       We compute their generating function,  finding elegant expressions for them in terms of the Cayley tree generating function. In \cite{LRT-bijective-visit} we put on the glasses of a bijective combinatorialist and added Enrica Duchi and Pablo Puerto to our team, and revisited our results. 
       In doing so, we establish connections between records in trees, endofunctions of finite sets and the \emph{genesis sequence}, the very first sequence in the Online Encyclopedia of Integer Sequences.
       \vspace{.3cm}

\noindent
{\bf Act III: }
After their exhausting journeys, the cars arrive one by one at their destination: a long street lined with plenty of parking spots. 
Since overtaking is impossible, at each step the fastest unblocked car reaches the street first. Too tired to be picky, cars simply park at the first available spot. 
\vspace{.15cm}

The main idea of Section \ref{se:parking} is that the parking scenario we just described provides an alternative framework to study parking functions, a perspective that leads to a record-preserving bijection between parking functions and rooted trees. Section \ref{se:parking} splits into three subsections.
In Section~\ref{Sec:WearyParking}, we formalize the parking scenario described in  Act III, where cars arrive at the parking lot sequentially and park in the first available spot. We refer to this process as \emph{weary parking}. 

The weary parking algorithm can be defined in terms of the \emph{priority-first search} algorithm. Here one starts with $\circ$ as the only visited node, and at each step visits the (previously unvisited) node with the smallest label that is adjacent to some previously visited node. 
This traversal is closely related to \emph{Prim's algorithm} for edge-weighted graphs, which we elaborate on in Section \ref{sec:final-comments}.  See also \cite[Chapter 4.3]{Sedgewick} for more details.

 How can we associate a parking function to a given weary parking process? Since drivers are tired and eager to park quickly, one might be tempted to say that all prefer spot 1. However, this is a poor choice, as we lose all the information about the weary parking process. On the other hand, since overtaking is impossible, the best drivers can hope for is to park immediately after the car in front. So we make this spot their preference. The bright side of this rule is that everyone is guaranteed a spot, and this preference function  always results in a parking function.

 Let $\mathcal{T}_n$ be the set of rooted trees of labelled at $[n] \cup \{\circ\},$ and rooted at $\circ$, and let $\mathcal{PF}_n$ be the set of parking functions of length $n$. 
 The key construction of this work is the \emph{weary bijection}, a bijective map from $\mathcal{T}_n$ to $\mathcal{PF}_n$ that sends a rooted tree $T$ to a parking function obtained from the weary parking process. To show that it is a bijection, we describe its inverse map. We also show that the weary permutation preserves records, as well as a host of other interesting statistics that we will describe later.
 
  The permutation describing the order in which the  cars end up parked under a parking  process is referred to as its \emph{bird's eye permutation}. 
  An important feature of the weary bijection is that, if $\pi_T$ is the parking function obtained as the image under the weary bijection of a rooted tree $T$, then, the bird’s eye permutations of $\pi_T$ and the one of the weary parking process of $T$  coincide. 
  Finally we remark that the bird's eye permutation is the inverse of the well-studied outcome permutation.

Some remarks on the history of parking functions. Konheim and Weiss  introduced them in their study of linear hashing, and showed, using generating functions, that there are \((n+1)^{n-1}\) parking functions of length \(n\), the same as the number of rooted trees with $n+1$ nodes, see  \cite{KW}. An unpublished argument by Henry O. Pollak, recounted in \cite{Foata-Riordan}, beautifully explains this result. 
 
  Several interesting bijections between parking functions and rooted trees are known. Many of them rely either on the recursive structure of trees and forests, Prüffer codes, or on the DFS (depth-first search) or the BFS (breath-first search) algorithms. These bijections are particularly powerful for  enumerative questions, including refinements of Cayley's formula by different statistics. They also provide efficient algorithms for generation, sorting, queue simulations, and data structure analysis.
   There is another family of bijections that involve intermediate Catalan objects, like Dyck paths or noncrossing partitions. These Catalan based bijections appear naturally in the study of algebraic or poset-theoretic questions, such as those arising in the study of diagonal harmonics, Weyl groups, or lattice theory. The reader may consult the excellent surveys \cite{  Stanley_parking,  Yan_handbook} and the references therein for more information.

The weary bijection introduced in this work belongs to the first family of bijections and explores another way of traversing a tree based on the priority search algorithm. 
Consider the weary parking process defined by a rooted tree $T$. Since overtaking is impossible, relabelling each car with the label of its final position always results in an increasing tree,  the \emph{priority tree} of $T$. The labels of the priority tree indicate the parking order of the cars under  weary parking (and not speed ranks). Moreover,  the $i^{th}$ car to park (represented by node $i$) prefers to park at the position immediately after the spot where its parent node parks. This observation is interpreted in Lemma \ref{le:action_weary_permutation} as saying that if the bird's eye permutation of $T$ is applied to the labels of its priority tree, we recover the original rooted tree $T$.

 Bijections obtained from different ways of transversing a tree keep track of different tree statistics:
Bijections based on BFS keep nodes of the same height close to one another. In contrast, parent-child pairs may appear far apart in the traversal. Thus, BFS favor proximity, and  can be used to keep track of the number of children of the root, \cite{Gessel-Seo}.  Bijections based on DFS, by contrast, preserve parent-child relationships, but lose information about the height of a node. Among other things, DFS has been used to keep track of tree inversions, \cite{Gessel-Sagan}.

Based on the priority search algorithm, the weary bijection preserves other interesting statistics.
In  Section~\ref{se:equidistribution},  we show that the weary bijection in record-preserving. 
 We build on the notion of priority tree to introduce some novel  statistics for rooted trees and parking function.  Then, we use the weary bijection to give a combinatorial proof of the join equidistribution of our statistics, and some other interesting statistics present in the literature. This is the content of Theorem \ref{thm:statistics}. 
 
 We also remark the existence of some well-studied pairs of statistics, including Kreweras's pair consisting of the inversion enumerator on parking functions and  the total displacement enumerator of rooted trees \cite{kreweras1980famille}, and  those studied in \cite{Stanley_parking} that are not equidistributed with the record statistic. This is elaborated in Remark \ref{re:kreweras_enumerator}, and the final comments.

  The equidistribution result also holds for restricted classes of trees/parking functions. This observation suggests the existence of a \emph{record duality} phenomenon for which we provide some interesting examples including  increasing trees and subexceedant functions. 

Observe that different trees can have the same priority tree. Thus, priority trees  induce a new partition of $\mathcal{PF}_n$, the \emph{priority partition}. Two parking functions are in the same class if they have the same priority tree.   

For example, on a one-way street (path graph), cars park according to the  order in which they enter the street, irrespective of speed ranks. Therefore, all $n!$ possibilities for the labels correspond to the same priority tree. 

 The set $\mathcal{PF}_n$ is often partitioned according to the sorted list of preferences. 
 Thus, the equivalence classes are labelled by Catalan objects, like Dyck paths or noncrossing partitions.
In Section \ref{weary_to_preference}, we show that the priority partition refines the Catalan partition. 

We close this work  with some enumerative results. Given a permutation $\omega \in \mathfrak{S}_n$, how many parking functions have $\omega$ as their bird’s-eye permutation? Given an increasing tree $T$, how many permutations have $T$ as their priority tree?  These questions are answered in Section \ref{se:hook}.  Note that the second question  asks for the sizes of the equivalence classes of the increasing tree partition of $\mathcal{PF}_n$.

We conclude this introduction by highlighting another elegant feature of the weary bijection: its inverse is equally natural and preserves the familiar “cars-that-want-to-park” intuition underlying parking functions. In contrast, most bijections in this area are markedly simpler in one direction than in the other.


\section{Preliminaries.} \label{se:blob}\label{se:record-blob}\label{suse:record}
We write  $[n]$ for  $\{1, 2, \ldots, n\}$ and $[n]_0$ for $\{0,  1, 2, \ldots, n\}$, both canonically ordered. As is customary, we sometimes use one-line notation to write permutations and parking functions, while at other times we use two-line notation, slightly abusing terminology.

Parking function are usually described by  a colourful scenario reminiscent of the one presented in the introduction.
Picture a one-way street with $n$ parking spaces, together with $n$ cars trying to park there. For each $i$, car $i$ has a preferred parking space $a_i$, $1 \le a_i \le n$.
Cars arrive on the street one after the other. When car $i$ enters the street, it first attempts to park in its preferred space $a_i$. If it is vacant, the car parks there. Otherwise, it continues down the street and attempts to park in the first available space it encounters. If all $n$ cars find a parking spot, the sequence of preferred spaces $\pi = (a_1, \ldots, a_n)$ is referred to as a \emph{parking function} of length $n$.

\begin{figure}[h!]
    \centering
    \resizebox{0.85\textwidth}{!}{
    \newcommand{\mzncar}[2][(0,0)]{%
\begin{scope}[shift={#1}]%
\draw[draw=black, fill=black, rounded corners=1.2ex, very thick] (0,.5) -- ++(0,1) -- ++(1,0.3) -- ++(3,0) -- ++(1.1,0.05) -- ++(0,-1.3) -- cycle;%
\draw[very thick, rounded corners=0.5ex, fill=black, thick] (1,1.74) -- ++(1,0.7) -- ++(1.6,0) -- ++(0.9,-0.7) -- cycle;%
\draw[fill=white, draw=white, thick] (1.25,.5) circle[radius=.55];%
\draw[fill=white, draw=white, thick] (4,.5) circle[radius=.55];%
\draw[fill=black,semithick] (1.25,.5) circle[radius=.35];%
\draw[fill=black,semithick] (4,.5) circle[radius=.35];%
\node[white] (label) at (2.7,1.5) {\bf \Huge #2};%
\end{scope}%
}%
\newlength{\xlen}%
\setlength{\xlen}{5.7cm}%
\newlength{\ylen}%
\setlength{\ylen}{8cm}%
\newcommand{\mzncars}[2]{%
    \foreach \label [count=\n from 0] in {#1} {%
        \ifdefstring{\label}{-}{%
        \phantom{\mzncar[(\n * \xlen, #2 * \ylen)]{\label}}%
        }{%
        \mzncar[(\n * \xlen - 0.3cm, #2 * \ylen)]{\label}%
        }%
        \coordinate (coo) at (\n * \xlen, #2*\ylen);%
        \draw[USred, line width=1mm] (coo) -- +(0.8*\xlen, 0);%
    }%
}%
\begin{tikzpicture}
\node (left) at (0, 0) {
\begin{tikzpicture}[scale=.5, xscale=-0.85, yscale=0.9]
\mzncars{-, -, -, -, -, -, -, -, -}{0}
\mzncars{-, -, -, -, 1, -, -, 2, -}{-1}
\mzncars{-, -, -, 3, 1, -, 4, 2, -}{-2}
\mzncars{-, -, 6, 3, 1, -, 4, 2, 5}{-3}
\mzncars{-, 8, 6, 3, 1, 7, 4, 2, 5}{-4}
\end{tikzpicture}
};
\qquad
\node (right) at (25, -2) {
\begin{tikzpicture}[scale=.5, xscale=-0.85, yscale=0.9]
\mzncars{-, -, -, -, 1, -, -, -, -}{-0}
\mzncars{-, -, -, 3, 1, -, -, 2, -}{-1}
\mzncars{-, -, -, 3, 1, -, 4, 2, 5}{-2}
\mzncars{-, -, 6, 3, 1, 7, 4, 2, 5}{-3}
\mzncars{9, 8, 6, 3, 1, 7, 4, 2, 5}{-4}
\end{tikzpicture}
};
\end{tikzpicture}
    }
    \caption{The parking function $ \pfunction{\pi = 5, 2, 5, 3, 1, 6, 1, 2, 6}$.}
    \label{fig:classical_parking}
\end{figure}

\begin{ex}
\label{ex:classical_parking_ex} The parking function 
 $\pi= \pfunction{5, 2, 5, 3, 1, 6, 1, 2, 6}$ has length 9.  The parking process defined by $\pi$ starts with the first car successfully parking at its preferred space $\pi(1) = 5$. The second car also gets lucky, and parks at its preferred space $\pi(2) = 2$. The third car, seeking the same location as the first one, discovers it claimed and advances to the next vacant space: spot 6. This process repeats until all cars have entered the street and parked as depicted in Figure \ref{fig:classical_parking}. 
  
\end{ex}

Rooted labelled trees of finite order are the other protagonist of this work. We assume that all trees are finite and labelled, and we omit these two adjectives. 
Unless stated otherwise, we assume that $T$ stands for a rooted tree with $n+1$ nodes, that the label of the root is $\circ$, and that  the set of labels for the remaining nodes is $[n]$. Moreover, we refer to $T$ as a tree of \emph{order} $n$, as we do not consider $\circ$ to be a real node. Alternatively, we say that $T$ has been \emph{labelled with $[n]_0$} and identify $\circ$ with 0.

A \emph{record of a rooted tree} is a node that has the largest label along the unique path from itself to the root.  Denote the set of records of \( T \) by \emph{\(\Rec(T)\)}, and its cardinality by \emph{\(\rec(T)\)}. Just as we do not consider $\circ$ to be a real node, when $T$ is rooted at $\circ$, we do not consider it to be a record either. Thus, $\Rec(T) \subseteq [n].$

The notion of record suggests a natural decomposition of a rooted tree in terms of simpler trees with their roots as the only record. A \emph{bonsai tree} is a rooted tree, labelled with positive integers, such that  its root is its only record.  In particular, this implies that node $\circ$ does not appear in any bonsai.
There is a natural way of obtaining a sequence of bonsai trees from a rooted tree \( T \). Delete from $T$ all edges joining a record node with its parent node. The result is a forest of bonsai trees.   The sequence obtained by  ordering the set of bonsais of $T$ increasingly according to the label of their roots will be referred to as the \emph{bonsai sequence} of $T$. 

\begin{figure}[h!]
\centering
\begin{subfigure}{0.25\textwidth}
    \centering
    \normalfont
\begin{tikzpicture}
    [
    level distance=0.6cm,
    level 1/.style={sibling distance=1.2cm},
    level 2/.style={sibling distance=0.75cm,},
    level 3/.style={sibling distance=0.75cm},
    cutedge/.style={dashed},
    every label/.style={inner sep=2pt,font=\tiny}
    ]
    \node[non-record, label={right:$\circ$}] (V-0) {}child {	
	    node[record, label={right:5}] (V-5) {}child {	
		    node[non-record, label={right:2}] (V-2) {}child {	
			    node[non-record, label={right:4}] (V-4) {}}}child {	
		    node[record, label={right:8}] (V-8) {}}}child {	
	    node[record, label={right:7}] (V-7) {}child {	
		    node[non-record, label={right:1}] (V-1) {}child {	
			    node[non-record, label={right:6}] (V-6) {}}child {	
			    node[record, label={right:9}] (V-9) {}}}child {	
		    node[non-record, label={right:3}] (V-3) {}}};
\end{tikzpicture}
    \caption{}
    \label{Fig:running_example}
\end{subfigure}
\begin{subfigure}{0.25\textwidth}
    \centering
    \begin{tikzpicture}
    [
    level distance=0.6cm,
    level 1/.style={sibling distance=1.2cm},
    level 2/.style={sibling distance=0.75cm,},
    level 3/.style={sibling distance=0.75cm},
    cutedge/.style={dashed},
    every label/.style={inner sep=2pt,font=\tiny}
    ]
    \node[non-record, label={right:$\circ$}] (V-0) {}
        child {node[record, label={right:5}] (V-5) {} 
            child {node[non-record, label={right:2}] (V-2) {}
                child {node[non-record, label={right:4}] (V-4) {}
                edge from parent[solid]}
            edge from parent[solid]}
            child {node[record, label={right:8}] (V-8) {}
            edge from parent[cutedge]}
        edge from parent[cutedge]}
        child {node[record, label={right:7}] (V-7) {}
            child {node[non-record, label={right:1}] (V-1) {}
                child {node[non-record, label={right:6}] (V-6) {}
                edge from parent[solid]}
                child {node[record, label={right:9}] (V-9) {}
                edge from parent[cutedge]}
            edge from parent[solid]}
            child {node[non-record, label={right:3}] (V-3) {}
            edge from parent[solid]}
        edge from parent[cutedge]};
\end{tikzpicture}
    \caption{}
    \label{Fig:record_edges_deletion}
\end{subfigure}
\begin{subfigure}{0.38\textwidth}
    \centering
    \begin{tikzpicture}[
    level distance=0.6cm,
    level 1/.style={sibling distance=0.75cm},
    level 2/.style={sibling distance=0.75cm},
]
\matrix[matrix of nodes, row sep=0.3cm, column sep=0.4cm, ampersand replacement=\&]{
    \node[record, label={right:5}] (V-5) {} 
        child {node[non-record, label={right:2}] (V-2) {}
            child {node[non-record, label={right:4}] (V-4) {}}}; \& 
    \node[record, label={right:7}] (V-7) {}
        child {node[non-record, label={right:1}] (V-1) {}
            child {node[non-record, label={right:6}] (V-6) {}}}
        child {node[non-record, label={right:3}] (V-1) {}}; \&
    \node[record, label={right:8}] (8) {};\& 
    \node[record, label={right:9}] (9) {};\\    
};
\end{tikzpicture}
    \caption{}
    \label{Fig:record_attachment_seq}
\end{subfigure}
\caption{(\subref*{Fig:running_example}) 
A rooted tree $T$ labelled with $[9]_0$  with record nodes \( \{5, 7, 8, 9\} \).
(\subref*{Fig:record_edges_deletion}) 
The tree $T$ with edges beginning in a record node deleted. 
(\subref*{Fig:record_attachment_seq}) The bonsai sequence of $T$. Its attachment sequence is $(\circ, 5, 1).$
}
\label{Fig:record_figures}
\end{figure}

 Different trees can share a bonsai sequence. To see this, observe that if we attempt to reconstruct a tree from a bonsai sequence $b=(B_1, B_2, \ldots, B_k)$, then the bonsai $B_1$ must be attached to $\circ$. However, bonsai $B_2$ could be attached to any of the nodes of $B_1$ or to $\circ$, and so on.  However, if in addition to the bonsai sequence of $T$, we register the parents of all the record nodes we have enough information to reconstruct $T$. 
We refer to the resulting sequence  as the \emph{attachment sequence} of $T$. 
An example of a rooted tree, together with its bonsai decomposition and attachment sequence is presented in Figure
\ref{Fig:record_figures}.
The \emph{record decomposition} of  $T$ is the ordered pair consisting of its bonsai sequence  and its attachment sequence. For more information, see \cite{Kreweras-Moszkowski, LRT-GF-records}.


\section{Weary Drivers Park Fast}
\label{se:parking}
\label{se:parking_weary}

\subsection{Weary parking}
\label{Sec:WearyParking}

Let $T$ be a rooted tree representing a bird's-eye view of the relative positions of cars, and imagine that the root $\circ$ of $T$ denotes the entry point to a street with plenty of parking space on one of its sides.
At each point of the process, the fastest unblocked car  reaches $\circ$ first, and parks at the nearest available spot.
We refer to this parking process as \emph{weary parking}. 

Observe that under weary parking, no car is totally unlucky. By construction, all cars end up finding a parking space. The bird's eye photo of the resulting sequence of parked cars describes a permutation in $\SS_n$ that we term  the \emph{bird's eye permutation} of $T$, $\omega_T$. Explicitly, $\omega_T(i)$ is the car that ends up parked at position $i$ after the parking process concludes. 
The definition of bird's eye permutation also applies  to any parking function $\pi$, in which case we write $\omega_\pi$.

During the weary parking process, cars arrive in clusters (the bonsais of $T$) led by the slowest vehicle of the cluster (the root of the bonsai). They are then linearized into a queue (park in a long street), headed by their record nodes. Moreover, since all ancestors of a record are smaller than it, the clusters (bonsais) arrive from earliest to latest in increasing order of their roots (record nodes).

\begin{ex}[Weary parking]
\label{ex:weary}
Let $T$ be the rooted tree depicted in Figure \ref{fig:enter-label}.  The weary parking process proceeds as follows.
The fastest unblocked car is 5, so it parks first. 
Then, car 2 becomes the fastest unblocked car and parks next. This, in turn, unblocks car 4, which then proceeds to park. 
Observe that the entire bonsai rooted at 5 parks first.

Car 7 becomes the fastest unblocked car. Thus, cars in its bonsai park sequentially: 7, then 1, then 3, then 6. 
 Only cars 8 and 9 remain unparked. Both are unblocked but 8 is faster, hence 8 parks first, leaving car 9 to be the last one to park.

\begin{figure}[h!]
    \centering
    \resizebox{\textwidth}{!}{\normalfont \begin{tikzpicture}[]
  \matrix[matrix of nodes, row sep=0.5cm, column sep=0.2cm] {
    \input{tikz.parking_arrival_1} & 
    \input{tikz.parking_arrival_2} & 
    \input{tikz.parking_arrival_3} & 
    \input{tikz.parking_arrival_4} &
    \input{tikz.parking_arrival_5} \\ 
    \input{tikz.parking_arrival_6} & 
    \input{tikz.parking_arrival_7} & 
    \input{tikz.parking_arrival_8} &
    \input{tikz.parking_arrival_9} & 
    \input{tikz.parking_arrival_10} \\ 
  };
\end{tikzpicture}}
    \caption{The weary parking process for the tree $T$ of our running example.}
    \label{fig:enter-label}
\end{figure}
The bird's eye permutation of $T$, written in one-line notation, is  
$ \omega_T =  \textcolor{CB.red}      {5 \ 2 \ 4 \ } 
   \textcolor{CB.orange}   {7\ 1 \ 3 \ 6\ } 
   \textcolor{CB.blue}     {8\ } 
   \textcolor{CB.green}    {9}.
        $
 Observe that a subword beginning at a record and ending just before the next one is the bird's eye permutation of the corresponding bonsai.
\end{ex}

Given two cars $i,j$ of $T$, one may ask how to tell who parks first. Clearly if $i$ is an ancestor of $j$ (or vice-versa), then $i$ parks before $j$ (and vice-versa). What if neither is an ancestor of the other?
\begin{lem} \label{lem:weary-order}
    Let $i,j$ be nodes of $T$ such that neither is the ancestor of the other. Let $k$ be the common ancestor of $i$ and $j$ at the largest height, and let $i',j'$ be the largest labels encountered on the paths from $i$ to (but not including) $k$ and $j$ to (but not including) $k$. Then $i$ parks before $j$ if and only if $i' < j'$.
    
\end{lem}

\begin{proof}
    Clearly neither $i$, nor $j$ are unblocked before $k$ is unblocked. After $k$ is unblocked, each of the ancestors of $i$ below $k$ have label at most $i'$, and thus park before $j'$. Therefore $i$ will be unblocked before $j$ parks, and will thus park before $j$.
\end{proof}

Let us note a few facts that follow immediately from this previous lemma. Firstly, we consider the case when the nodes $i,j$ are from the same bonsai. Then the common ancestor $k$ must also be in the same bonsai, and so we see that cars in the same bonsai park in the same order regardless of the rest of the tree $T$. Secondly, we consider the case when $i,j$ are in different bonsais with respective roots $r_1,r_2$ so that $r_1 < r_2$. Then it is not hard to see that $i' \leq r_1$ and $j' = r_2 > i'$, so that $i$ parks before $j$. 

The second of the above facts justifies that under the weary parking process that cars park one bonsai at a time, in increasing order of the roots of the bonsais. On the other hand, the first fact confirms that the parking order in each bonsai is independent of the others. Putting these two facts together we see that the bird's eye permutation does not depend on the attachment sequence, only on the bonsai sequence. 

By artificially rooting a given bonsai $B$ at $\circ$, one may consider the weary parking process restricted to the bonsai $B$. It follows immediately from Lemma \ref{lem:weary-order} that the vertices of $B$ must park in the same order in this restricted process as they would if $B$ was a bonsai of any rooted Cayley tree $T$  on $[n]_0$. 
Therefore, we may also describe the weary parking process recursively.

Define the \emph{bonsai word} of a rooted tree $T$ as the formal word that results from concatenating the bonsais appearing in the decomposition of \( T \), and increasingly ordering them according to the labels of  their roots.  The process of weary parking can be recursively described as follows. Let $W$ be the bonsai word of $T$. Then, first write the label of the root of each bonsai in this word immediately to its left, and relabel the root of the bonsai as $\circ$, and then discard each bonsai consisting of a single node, the root $\circ$.  Afterwards, apply this recursion to the remaining bonsai trees of the resulting word.  
Since $T$ is finite, at some point there will be no trees left in this expression, and the procedure stops. 
The result is the bird's eye permutation $\omega_T$ of $T.$ We summarize this observation in the following lemma. 

\begin{lem}
   The word obtained at the end of the recursive algorithm applied to $T$ is the bird's eye permutation $\omega_T$.
\end{lem} 

 \begin{ex} \label{ex:weary_recursive}
\newcommand{\treeI}{
    \begin{tikzpicture}[level distance=0.5cm, baseline=-0.5cm]
    \node[non-record, label={right:$\circ$}] (V-5) {} 
        child {node[record, label={right:2}] (V-2) {}
            child {node[record, label={right:4}] (V-4) {}}};
    \end{tikzpicture}
}

\newcommand{\treeII}{
    \begin{tikzpicture}[
        level distance=0.6cm,
        level 1/.style={sibling distance=0.75cm},
        level 2/.style={sibling distance=0.75cm}, 
        baseline=-0.5cm
        ]
    \node[non-record, label={right:$\circ$}] (V-7) {}
        child {node[record, label={right:1}] (V-1) {}
            child {node[record, label={right:6}] (V-6) {}}}
        child {node[record, label={right:3}] (V-1) {}}; 
    \end{tikzpicture}
}

\newcommand{\singlenode}{
\begin{tikzpicture}
\node[non-record, label={right:$\circ$}] (root) {};
\end{tikzpicture}
}

 Using the same tree as in Example \ref{ex:weary}, we perform the recursive algorithm, whose bonsai sequence appears in Figure \ref{Fig:record_attachment_seq}. After a first iteration, we obtain
\begin{center}
\normalfont
(\memph{5}, \, \treeI \hspace{-1ex}, \, \memph{7}, \, \treeII \hspace{-1ex}, \, \memph{8}, \, \memph{9})
\end{center}
(the trivial bonsais with roots 8 and 9 have been erased during part (2) of this iteration). Since all remaining trees are increasing,  all their bonsais consist of a single node. Thus, after the second step we obtain
\begin{center}
\normalfont
(5, \memph2, \memph4, 7, \memph1, \memph3, \memph6, 8, 9).
\end{center}
As there are no trees left in the word, this is its bird's eye permutation.
 \end{ex}

The recursive description  has some immediate consequences that we summarize in the following corollary. 

\begin{cor} \label{cor:concatenation}  Let $T$ be a rooted tree of order $n.$
Then, the bird's eye permutation $\omega_T$ of $T$ is the concatenation of the bird's eye permutations of the bonsais, respecting the order in which they appear in the bonsai sequence of $T.$ In particular, $\omega_T$ does not depend on the attachment sequence of $T$.
\end{cor}

The weary parking process can be described as follows:
Start at the root $\circ$, and at each step visit the node with the lowest label that is adjacent to one of the previously visited nodes. In terms of the \emph{priority first search} algorithm, the bird's eye permutation of $T$ is the tree transversal obtained by doing a priority first search on $T$ starting at $\circ$. When we want to emphasize this point, we refer to the bird's eye permutation of $T$ as its \emph{priority transversal}.

\subsection{The weary bijection. From rooted trees to parking functions}
\label{weary_to_preference}

Now that we know how cars in a given rooted tree park, we ask ourselves how to associate a parking function to this process. This amounts to assigning a parking preference to every car. 
Since drivers are tired and eager to park quickly, one might be tempted to say that all prefer spot 1. However, this choice makes it impossible to establish a bijection between trees and parking functions. On the other hand, since overtaking is impossible, the best drivers can hope for is to park immediately after the car in front. The bright side of this rule is that everyone is guaranteed a spot.

We formalize this observation as follows. The \emph{parent function} of a rooted tree \( T \), denoted by  \( p_T \), is  the function that sends a node of \( T \) other than the root to its parent node, and the root to itself. 
Let $\omega_T$ be the bird's eye permutation of the weary parking process defined by $T$.  
The \emph{preference} $\pi_T(i)$ of  car $i$  is defined as
\begin{align}
\label{Def:weary_preference}
\Pref(i) = \omega_T^{-1}(p_T(i)) + 1.
\end{align}
The reasoning behind this definition is this. First, car $i$ looks at the spot where its parent $p_T(i)$ parks. Since cars park consecutively starting at spot 1, the arrival order coincides with the label of the parking spots. Thus, the parent of car $i$ takes spot  $\omega_T^{-1}(p_T(i))$, and  car $i$ decides that it would like to park at the spot immediately after it, $ \omega_T^{-1}(p_T(i)) + 1$.

Define the \emph {priority vector} of a rooted tree $T$ as 
\begin{equation}
\label{eqn:pv}
\pv(\pi_T) = \pi_T \circ \omega_T.
\end{equation}
It is clear that $\pv(\pi)(i)$ is the spot where the car that ends up parked at position $i$ would have liked to park. Moreover, the car that ends up parked at position $i$ would have really liked to park  in position $\pv(\pi)(i).$

We say that a function $f : [n] \to [n]$ is \emph{subexceedant} if $f(i) \le i$ for all $i$ in $[n].$ Note that the priority vector  $\pv(\pi_T)$ of $\pi_T$ is a subexceedant function, and that a subexceedant function is always a parking function. As customary, the symmetric group acts on a function $f: [n] \to [n]$ by the rule $\omega \cdot f(i)= f\big( \omega^{-1}(i) \big)$ for any $\omega \in \mathbb S_n$.  
As the set of parking functions is closed under the action of the symmetric group, Equation \eqref{eqn:pv} implies that $\pi_T$ is a parking function.

\begin{ex} 
\label{ex:pv_canonical_form}
The priority vector of our running example,
$
\pi=
\pfunction{5,2,5,3,1,6,1,2,6}      
$ is $
\begin{matrix}
\pv(\pi) = \pfunction{1, 2, 3, 1, 5, 5, 6, 2, 6}.
\end{matrix}
 $ It carries the following information.
 After the parking process finishes, the cars parked at positions 1 and 4 wanted spot 1; cars at 2 and 8 wanted spot 2; car 3 wanted spot 3; and so on, up to car 9, which wanted spot 6.
    
\end{ex}

Cars that park in their preferred spot are known as
\emph{lucky cars}. After all, luck is wanting to do what we must do anyway. The fixed points of  $\pv(\pi)$ are the \emph{lucky spots} of $\pi$, that is, the positions occupied by the lucky cars of $\pi$.

We want to stress that   the sequence obtained from  $\pi $  (written using one-line notation) by rearranging it as a non-decreasing sequence is in general a different sequence than the priority vector $\pv(\pi)$. See Example \ref{ex:pv_canonical_form}.

\begin{thm}[The weary bijection] Let $\mathcal{T}_n$ be the set of rooted trees of order $n$ and $\mathcal{PF}_n$ be the set of parking functions of length $n$. Then, the  map $\mathcal{T}_n \to \mathcal{PF}_n $ such that  $ T \mapsto  \pi_T $,  where $\pi_T$ is the parking function defined by the weary parking process and  Eq.(\ref{Def:weary_preference}),
defines a bijection between $ \mathcal{PF}_n$ and $\mathcal{T}_n$. 
\end{thm}

 We will refer to this bijection as the \emph{weary bijection}. In the next section, we  show that the weary bijection is indeed bijective by constructing its inverse map.

\begin{ex}
\label{ex:running_as_image_of_tree} Let $T$ be the rooted tree of Figure \ref{fig:enter-label}. The parking function of $T$ (its image under the weary bijection)  is  
 $\pi_T= \pfunction{5, 2, 5, 3, 1, 6, 1, 2, 6}$, the parking function of Examples \ref{ex:classical_parking_ex} and \ref{ex:weary_recursive}.  The bird's eye permutations of $\Pref$ and $T$ are both equal to $  
\pfunction{5, 2, 4, 7, 1, 3, 6, 8, 9}
$. 
\end{ex}

\begin{cor}
The weary bijection restricts to a bijection between the set of increasing trees of order $n$, and the set of subexceedant functions of length $n.$
\end{cor}

A crucial property of the weary parking process on a tree $T$ is that it always yields the same parking configuration as the parking function $\pi_T$.

\begin{lem}
\label{lem: omega_T = omega_pi_T}
    For any rooted tree $T$, the bird's eye permutation of $T$ coincides with the bird's eye permutation of its preference sequence. That is, $\omega_T = \omega_{\pi_T}$.
\end{lem}
\begin{proof}
    We proceed by induction. For the base case, note that $\omega_T(\circ) = \omega_{\pi_T}(\circ)$ by convention. 
    Next, suppose that we have already determined that $\omega_T(i) = \omega_{\pi_T}(i)$ for all $i \in [k]_0$, denote $U = \{\omega_T(\circ), \omega_T(1), \dots, \omega_{T}(k)\}$ and let us show that $\omega_T(k+1) = \omega_{\pi_T}(k + 1)$.
    
    The element $\omega_T(k+1)$ is the node visited by the priority search of $T$ at step $k+1$, so it must be the minimum $i$ that is unblocked but not visited after step $k$. More rigorously put, $\omega_T(k+1)$ is the minimum element $i$ in $[n] \setminus U$ such that $p_T(i) \in U$. 

    On the other hand, the element $\omega_{\pi_T}(k+1)$ is the car that parks at space $k + 1$ for the parking function $\pi_T$. Naturally, $\omega_{\pi_T}(k+1)$ must lie in the set $[n] - U$, since the cars in $U$ park in the first $k$ spots. Furthermore, in order for car $\omega_{\pi_T}(k+1)$ to find space $k+1$ unclaimed, it must be the lowest element $i$ in $[n] \setminus U$ such that $\pi_T(i) \leq k+1$. By definition of $\pi_T$, this last inequality is equivalent to $\omega_T^{-1}(p_T(i)) \leq k$, and henceforth $\omega_{\pi_T}(k + 1)$ is the minimum $i$ in $[n] \setminus U$ such that $p_{T}(i) \in U$. This shows $\omega_T(k+1) = \omega_{\pi_T}(k + 1)$, as claimed.
\end{proof}

The results of the past two sections (Corollary \ref{cor:concatenation} in conjunction with Lemma \ref{lem: omega_T = omega_pi_T}) imply also that the weary bijection is record preserving. We offer a more direct proof of this in Section \ref{sec:record_bijection}.

\subsection{The inverse weary bijection. From parking functions to trees.}

In this section, we describe explicitly the \emph{inverse weary bijection},  the inverse of the weary bijection. Note that it is a remarkably simple procedure, especially when compared with other known bijections between rooted trees and parking functions.

 Let $\pi$ be a parking function of length $n$, and let $\omega_\pi$ be its bird's eye permutation, a permutation of $\SS_n$. 
 The inverse weary bijection maps $\pi$ to a rooted  tree $T$ of order $n$, that we call  the \emph{parking tree of $\pi$} and denote by \emph{$\PT(\pi)$}. It it is constructed as follows. Let  $S_k$ be the set of cars with preference $k$ in $\pi$. The elements of $S_k$ will become siblings (share a parent) on $T$.  The tree $T$ is constructed by attaching  the nodes of $S_1$  to the root, and the nodes of $S_k$ to $\omega_{k-1}$ (the  car parked at the position $k-1$) when $k>1$.

This procedure defines the inverse of the weary bijection. To see this, note that under the weary bijection, a set of siblings of a rooted tree 
 $T$ will be assigned the same preference, $\omega_T^{-1}(p_T(i))+1$, which is the position where $p_T(i)$ parks.
 
Note that any rooted tree is then the parking tree of a unique parking function. Formally, $\PT(\pi)$ is defined by the parent function $p_{\pi} : [n]_0 \to [n]_0$ given by $p_{\pi}(i) = \circ$ if $i = \circ$ or $\pi_i = 1$, and $p_{\pi}(i) = \omega_{\pi}(\pi_i-1)$ otherwise. 

The inverse of the weary bijection can be read from the \emph{parking blueprint} of $\pi$, a diagram  constructed as follows. First, we draw $\len(\pi)$ parking spots, and then, reading $\pi$ from left-to-right, we assign a car to each spot following the instruction contained in the parking function $\pi$: if position $\pi_k$ is empty, we draw car $k$ at $\pi_k$. Otherwise, car $k$ parks at  $\omega_\pi(i)$, the first empty position it encounters, and we draw car $k$ there. In addition, we draw an arc from $\pi_k$ to $\omega_\pi(i)$, taking care that the different arcs do not intersect. See Figure \ref{subfig:blueprint} for the parking blueprint of our running example.

The inverse parking bijection can be computed using parking blueprints.
To construct the tree that corresponds to $\pi$, first draw  the parking blueprint of $\pi$, and identify nodes with  cars, and observe that the set $S_k$ of cars with preference $k$ consist of all cars whose arc starts at position $k$. The nodes of $S_k$ are then attached to the root, when $k=1$, or to the car parked at position $k-1$, if $k>1$. This is illustrated in Figure \ref{subfig:blueprint to tree}.

\begin{figure}[h!]
    \centering
    \begin{subfigure}{0.4\textwidth}
        \centering
        \resizebox{\textwidth}{!}{
            \newcommand{\mzncar}[2][(0,0)]{%
\begin{scope}[shift={#1}]%
\draw[draw=black, fill=black, rounded corners=1.2ex, very thick] (0,.5) -- ++(0,1) -- ++(1,0.3) -- ++(3,0) -- ++(1.1,0.05) -- ++(0,-1.3) -- cycle;%
\draw[very thick, rounded corners=0.5ex, fill=black, thick] (1,1.74) -- ++(1,0.7) -- ++(1.6,0) -- ++(0.9,-0.7) -- cycle;%
\draw[fill=white, draw=white, thick] (1.25,.5) circle[radius=.55];%
\draw[fill=white, draw=white, thick] (4,.5) circle[radius=.55];%
\draw[fill=black,semithick] (1.25,.5) circle[radius=.35];%
\draw[fill=black,semithick] (4,.5) circle[radius=.35];%
\node[white] (label-#2) at (2.7,1.5) {\bf \Huge #2};%
\end{scope}%
}%
\setlength{\xlen}{5.7cm}%
\setlength{\ylen}{8cm}%
\newcommand{\mzncars}[2]{%
    \foreach \label [count=\n from 0] in {#1} {%
        \ifdefstring{\label}{-}{%
        \phantom{\mzncar[(\n * \xlen, #2 * \ylen)]{\label}}%
        }{%
        \mzncar[(\n * \xlen - 0.3cm, #2 * \ylen)]{\label}%
        }%
        \coordinate (coo-\label) at (\n * \xlen, #2*\ylen);%
        \draw[USred, line width=1mm] (coo-\label) -- +(0.8*\xlen, 0);%
    }%
}%

\tikzset{
    bp-circ/.style={
        circle,
        draw=USred,
        fill=USred,
        inner sep=3pt,
        label={#1},
        solid
    },
    bp-line/.style={
        line width=2.5pt, rounded corners=25pt, USred
    },
  every label/.style={font=\fontsize{28}{48}\selectfont}
}

\begin{tikzpicture}
\normalfont
\def\ysep{3}
\node (right) at (0, 0) {
\begin{tikzpicture}[scale=.5, xscale=-0.85, yscale=0.9]
\mzncars{9, 8, 6, 3, 1, 7, 4, 2, 5}{0}

\coordinate (mid-7) at ($(label-7)+(0, 1.5*\ysep)$);
\coordinate[very thick] (end-7) at ($(label-5) + (0, 1.5*\ysep)$);
\draw[bp-line]  (label-7) -- (mid-7) -- (end-7);
\node[bp-circ, label={left:7}] at (end-7) {};
\coordinate (mid-3) at ($(label-3)+(0, 2.5*\ysep)$);
\coordinate[very thick] (end-3) at ($(label-1) + (0, 2.5*\ysep)$);
\draw[bp-line]  (label-3) -- (mid-3) -- (end-3);
\node[bp-circ, label={left:3}] at (end-3) {};
\coordinate (mid-6) at ($(label-6)+(0, 3.5*\ysep)$);
\coordinate[very thick] (end-6) at ($(label-3) + (0, 3.5*\ysep)$);
\draw[bp-line]  (label-6) -- (mid-6) -- (end-6);
\node[bp-circ, label={left:6}] at (end-6) {};

\coordinate (mid-8) at ($(label-8)+(0, 4.5*\ysep)$);
\coordinate[very thick] (end-8) at ($(label-2) + (0, 4.5*\ysep)$);
\draw[bp-line]  (label-8) -- (mid-8) -- (end-8);
\node[bp-circ, label={left:8}] at (end-8) {};

\coordinate (mid-9) at ($(label-9)+(0, 5.5*\ysep)$);
\coordinate[very thick] (end-9) at ($(label-3) + (0, 5.5*\ysep)$);
\draw[bp-line]  (label-9) -- (mid-9) -- (end-9);
\node[bp-circ, label={left:9}] at (end-9) {};
\end{tikzpicture}
};
\end{tikzpicture}
        }
        \caption{}
    \label{subfig:blueprint}
    \end{subfigure}
    \qquad
    \begin{subfigure}{0.4\textwidth}
        \centering
        \normalfont
\begin{tikzpicture}
    [
    level distance=0.6cm,
    level 1/.style={sibling distance=1.2cm},
    level 2/.style={sibling distance=0.75cm,},
    level 3/.style={sibling distance=0.75cm},
    cutedge/.style={dashed},
    every label/.style={inner sep=2pt,font=\tiny}
    ]
    \node[non-record, label={right:$\circ$}] (V-0) {}child {	
	    node[record, label={right:5}] (V-5) {}child {	
		    node[non-record, label={right:2}] (V-2) {}child {	
			    node[non-record, label={right:4}] (V-4) {}}}child {	
		    node[record, label={right:8}] (V-8) {}}}child {	
	    node[record, label={right:7}] (V-7) {}child {	
		    node[non-record, label={right:1}] (V-1) {}child {	
			    node[non-record, label={right:6}] (V-6) {}}child {	
			    node[record, label={right:9}] (V-9) {}}}child {	
		    node[non-record, label={right:3}] (V-3) {}}};
\end{tikzpicture}
        \caption{}
    \label{subfig:blueprint to tree}
    \end{subfigure}
    \caption{(a) The parking blueprint of the parking function $\pi =
\pfunction{5,2,5,3,1,6,1,2,6}      
$ of Examples \ref{ex:classical_parking_ex} and 
\ref{ex:running_as_image_of_tree}. (b) The tree corresponding to $\pi$ under the inverse weary bijection.}
    \label{fig:blueprint}
\end{figure}

\begin{proof}[Proof that the weary map is a bijection]
We show that the maps 
\[
T \mapsto \pi_T \qquad \text{and} \qquad \pi \mapsto \PT(\pi)
\]
are inverses of each other. 
Recall that these maps are defined by
\[
\pi_T(i) = \omega_T^{-1}(p_T(i)) + 1 \qquad \text{and} \qquad p_\pi(i) = \omega_\pi(\pi(i) - 1),
\]
respectively, 
where we are writting $\omega_T(\circ) = \pi(\circ) = \omega_\pi(\circ) = \circ$. 
The result follows noting that $\omega_T = \omega_{\pi_T}$ (Lemma \ref{lem: omega_T = omega_pi_T}). 

\end{proof}

To summarize, observe that the parking function $\pi_T = (a_1,\dots,a_n)$ associated to the rooted tree $T$  can be computed concurrently with the bird's eye permutation: after a car $j$ is parked in some spot $k$, one can immediately assign the preferences of all of the children of $j$ to be $k+1$. 
The recursive process therefore may be summarized simply as: 
``park the unblocked car with the smallest label in the next available space $k$, unblock all of its children, set their preference to $k+1$, and repeat until all cars have been parked.'' The recursion begins with all children of $\circ$ unblocked with preference equal to $1$.
The reverse process may also be summarized by simply describing the parent function as: ``the parent of node (car) $i$ with preference $a_i$ is the node (car) $j$ that parked in spot $a_i - 1$''.

\section{The priority tree of a parking function.}
\label{sec:pf2tree}
 \label{se_arrival_tree}

Let  $\pi$ be a parking function and let  $\pv(\pi) = \sigma_1\sigma_2 \dots \sigma_{n}$ be the priority vector of $\pi$, and let  $p_T : [n]_0 \to [n]_0$ be its parent map. 
The \emph{priority tree} of $\pi$,  \emph{$\CT(\pi)$}, 
is the increasing tree that has as parent function the subexceedant function:
\begin{equation}
\label{eqn:priority tree map}
    p_T(i) = \sigma_i - 1. 
\end{equation} 

In the weary parking process defined by a rooted tree \(T\), overtaking is impossible. The effect of Equation~\ref{eqn:priority tree map} is to relabel each car by its final parking position, yielding an increasing tree.
 The position of node $i$ in $\CT(\pi)$ tells us \emph{when} car $\omega_\pi(i)$ parked under the weary parking process defined by $\PT(\pi)$.  

\begin{lem}
\label{le:nodes_parking_tree}
In the weary parking process on $\CT(\pi)$, car $i$ is the $i^{\text{th}}$ to park. Moreover, the $i^{th}$ car to park prefers to occupy the parking spot immediately following the one in which the $p_T(i)^{th}$ car parked.

\end{lem}

\begin{proof}
We proceed by induction on $i$. Since $\pv$ is a subexceedent function, $p_T(1) = \sigma_1 - 1 \leq 0$, so $p_T(1) = \circ$. Therefore, at the start of the weary parking process, the car $1$ is the fastest unblocked car and thus parks first. Now assume that the result holds for each $i=1,\dots,k-1$ for some integer $k \geq 2.$ Then, $p_T(k) = \sigma_{k} - 1 \leq k-1$. Since cars $1,\dots,k-1$ have already parked in the first $k-1$ rounds (by the inductive hypothesis), and since the parent of $k$ is in $\{1,\dots,k-1\}$, car $k$ is the fastest unblocked car in round $k$, and is thus the $k$th car to park, as required. The second statement is immediate.
\end{proof}

\begin{ex}[The parking tree of $\pi$]
\label{ex_parkingtree}
Let $
\pi=
\pfunction{5, 2, 5, 3, 1, 6, 1, 2, 6}
$ be the parking  function of  Example \ref{ex:classical_parking_ex}.
The parking tree of $\pi$ is then obtained by  making the bird's eye permutation $\omega_\pi$ (computed in Example \ref{ex:classical_parking_ex}) act on the nodes of $\CT(\pi)$. The result appears in Figure \ref{fig:parking_tree}.
\begin{figure}[h!]
    \centering
    \hspace*{\fill}%
    \begin{subfigure}[b]{0.35\textwidth}
        \centering
        \normalfont
\begin{tikzpicture}
    [
    level distance=0.6cm,
    level 1/.style={sibling distance=1.2cm},
    level 2/.style={sibling distance=0.75cm,},
    level 3/.style={sibling distance=0.75cm},
    cutedge/.style={dashed},
    every label/.style={inner sep=2pt,font=\tiny}
    ]
    \node[record, label={right:$\circ$}] (V-0) {}child {	
	    node[record, label={right:1}] (V-5) {}child {	
		    node[record, label={right:2}] (V-2) {}child {	
			    node[record, label={right:3}] (V-4) {}}}child {	
		    node[record, label={right:8}] (V-8) {}}}child {	
	    node[record, label={right:4}] (V-7) {}child {	
		    node[record, label={right:5}] (V-1) {}child {	
			    node[record, label={right:7}] (V-6) {}}child {	
			    node[record, label={right:9}] (V-9) {}}}child {	
		    node[record, label={right:6}] (V-3) {}}};
\end{tikzpicture}
        \caption{$\CT(\pi)$}
        \label{fig:arrival_tree}
    \end{subfigure}
    \begin{subfigure}[b]{0.35\textwidth}
        \centering
        \normalfont
\begin{tikzpicture}
    [
    level distance=0.6cm,
    level 1/.style={sibling distance=1.2cm},
    level 2/.style={sibling distance=0.75cm,},
    level 3/.style={sibling distance=0.75cm},
    cutedge/.style={dashed},
    every label/.style={inner sep=2pt,font=\tiny}
    ]
    \node[non-record, label={right:$\circ$}] (V-0) {}child {	
	    node[record, label={right:5}] (V-5) {}child {	
		    node[non-record, label={right:2}] (V-2) {}child {	
			    node[non-record, label={right:4}] (V-4) {}}}child {	
		    node[record, label={right:8}] (V-8) {}}}child {	
	    node[record, label={right:7}] (V-7) {}child {	
		    node[non-record, label={right:1}] (V-1) {}child {	
			    node[non-record, label={right:6}] (V-6) {}}child {	
			    node[record, label={right:9}] (V-9) {}}}child {	
		    node[non-record, label={right:3}] (V-3) {}}};
\end{tikzpicture}
        \caption{$\PT(\pi)$}
        \label{fig:parking_tree}
    \end{subfigure}
    \hspace*{\fill}%
    \caption{The priority and parking trees of the parking function $\pi$ of our running example. }
    \label{fig:trees_combined}
\end{figure}
\end{ex}

It is convenient to endow our rooted trees with a \emph{plane tree} structure where the children of each node are ordered increasingly. Then, the following relation between the priority tree and the parking tree of a parking function comes to light.

\begin{lem}
\label{le:action_weary_permutation}\label{se_parking_tree}
Let $\pi$ be a parking function $\pi$ with bird's eye permutation $\omega_\pi$ of $\pi$. Then, 
\begin{align}
\label{eqn:parking_tree_def}
\PT(\pi) = \omega_{\pi}(\CT(\pi)).
\end{align}
where $\omega_{\pi}$  acts naturally on $\CT(\pi)$ by permuting its labels. Moreover, the plane structure of $\PT(\pi)$, and the one inherited from $\CT(\pi)$ by this action coincide.
\end{lem}

\begin{proof}
After Lemma \ref{le:nodes_parking_tree}, it only remains to be shown that  the plane structure of the parking tree of $\pi$ is the same as the one that it inherits from the plane structure of the priority tree: if the ordered children of a node $i$ in $\CT(\pi)$ are $x_1, x_2, \dots, x_k$, then the ordered children of $\omega_\pi(i)$ in $\PT(\pi)$ are $\omega_\pi(x_1), \omega_\pi(x_2), \dots, \omega_\pi(x_k)$. 
Both sequences of nodes  $x_1,$ $x_2,$ $\dots,$ $x_k$ and $\omega_\pi(x_1),$ $\omega_\pi(x_2),$ $\dots,$ $\omega_\pi(x_k)$, are increasing. The first, as a consequence of the definition of the plane structure of $\CT(\pi)$. The second, because we are acting with the bird's eye permutation: Since $\omega_\pi(x_1),$ $\omega_\pi(x_2),$ $\dots,$ $\omega_\pi(x_k)$ share a parent, by construction, they prefer the same spot, hence they have to attempt to park in the order implied by their labels to result in the bird's eye permutation $\omega_\pi$. 
\end{proof}

Lemma \ref{le:action_weary_permutation} allows us to define the \emph{priority tree} of a rooted tree $T$ as $
  \pt(T) = \omega_T^{-1}(T)
$. 

Define the \emph{priority blueprint} of $\pi$ as the blueprint obtained  by relabeling each car by its parking position. However, these labels are redundant as they are always $1, 2, \ldots, n$
written from left to right. 
A convenient way of computing the inverse parking bijection is to start by computing its priority blueprint of $\pi$,  reconstructing from its priority tree from it, and then use Lemma \ref{le:action_weary_permutation}.

\begin{ex}
The parking blueprint of our running example and its priority blueprint.
\begin{figure}[h!]
    \centering
    \begin{subfigure}{0.4\textwidth}
        \centering
        \resizebox{\textwidth}{!}{
            \newcommand{\mzncar}[2][(0,0)]{%
\begin{scope}[shift={#1}]%
\draw[draw=black, fill=black, rounded corners=1.2ex, very thick] (0,.5) -- ++(0,1) -- ++(1,0.3) -- ++(3,0) -- ++(1.1,0.05) -- ++(0,-1.3) -- cycle;%
\draw[very thick, rounded corners=0.5ex, fill=black, thick] (1,1.74) -- ++(1,0.7) -- ++(1.6,0) -- ++(0.9,-0.7) -- cycle;%
\draw[fill=white, draw=white, thick] (1.25,.5) circle[radius=.55];%
\draw[fill=white, draw=white, thick] (4,.5) circle[radius=.55];%
\draw[fill=black,semithick] (1.25,.5) circle[radius=.35];%
\draw[fill=black,semithick] (4,.5) circle[radius=.35];%
\node[white] (label-#2) at (2.7,1.5) {\bf \Huge #2};%
\end{scope}%
}%
\setlength{\xlen}{5.7cm}%
\setlength{\ylen}{8cm}%
\newcommand{\mzncars}[2]{%
    \foreach \label [count=\n from 0] in {#1} {%
        \ifdefstring{\label}{-}{%
        \phantom{\mzncar[(\n * \xlen, #2 * \ylen)]{\label}}%
        }{%
        \mzncar[(\n * \xlen - 0.3cm, #2 * \ylen)]{\label}%
        }%
        \coordinate (coo-\label) at (\n * \xlen, #2*\ylen);%
        \draw[USred, line width=1mm] (coo-\label) -- +(0.8*\xlen, 0);%
    }%
}%

\tikzset{
    bp-circ/.style={
        circle,
        draw=USred,
        fill=USred,
        inner sep=3pt,
        label={#1},
        solid
    },
    bp-line/.style={
        line width=2.5pt, rounded corners=25pt, USred
    },
  every label/.style={font=\fontsize{28}{48}\selectfont}
}

\begin{tikzpicture}
\normalfont
\def\ysep{3}
\node (right) at (0, 0) {
\begin{tikzpicture}[scale=.5, xscale=-0.85, yscale=0.9]
\mzncars{9, 8, 6, 3, 1, 7, 4, 2, 5}{0}

\coordinate (mid-7) at ($(label-7)+(0, 1.5*\ysep)$);
\coordinate[very thick] (end-7) at ($(label-5) + (0, 1.5*\ysep)$);
\draw[bp-line]  (label-7) -- (mid-7) -- (end-7);
\node[bp-circ, label={left:7}] at (end-7) {};
\coordinate (mid-3) at ($(label-3)+(0, 2.5*\ysep)$);
\coordinate[very thick] (end-3) at ($(label-1) + (0, 2.5*\ysep)$);
\draw[bp-line]  (label-3) -- (mid-3) -- (end-3);
\node[bp-circ, label={left:3}] at (end-3) {};
\coordinate (mid-6) at ($(label-6)+(0, 3.5*\ysep)$);
\coordinate[very thick] (end-6) at ($(label-3) + (0, 3.5*\ysep)$);
\draw[bp-line]  (label-6) -- (mid-6) -- (end-6);
\node[bp-circ, label={left:6}] at (end-6) {};

\coordinate (mid-8) at ($(label-8)+(0, 4.5*\ysep)$);
\coordinate[very thick] (end-8) at ($(label-2) + (0, 4.5*\ysep)$);
\draw[bp-line]  (label-8) -- (mid-8) -- (end-8);
\node[bp-circ, label={left:8}] at (end-8) {};

\coordinate (mid-9) at ($(label-9)+(0, 5.5*\ysep)$);
\coordinate[very thick] (end-9) at ($(label-3) + (0, 5.5*\ysep)$);
\draw[bp-line]  (label-9) -- (mid-9) -- (end-9);
\node[bp-circ, label={left:9}] at (end-9) {};
\end{tikzpicture}
};
\end{tikzpicture}
        }
        \caption{}
    \end{subfigure}
    \qquad
    \begin{subfigure}{0.4\textwidth}
        \centering
        \resizebox{\textwidth}{!}{
            \newcommand{\mzncar}[2][(0,0)]{%
\begin{scope}[shift={#1}]%
\draw[draw=black, fill=black, rounded corners=1.2ex, very thick] (0,.5) -- ++(0,1) -- ++(1,0.3) -- ++(3,0) -- ++(1.1,0.05) -- ++(0,-1.3) -- cycle;%
\draw[very thick, rounded corners=0.5ex, fill=black, thick] (1,1.74) -- ++(1,0.7) -- ++(1.6,0) -- ++(0.9,-0.7) -- cycle;%
\draw[fill=white, draw=white, thick] (1.25,.5) circle[radius=.55];%
\draw[fill=white, draw=white, thick] (4,.5) circle[radius=.55];%
\draw[fill=black,semithick] (1.25,.5) circle[radius=.35];%
\draw[fill=black,semithick] (4,.5) circle[radius=.35];%
\node[white] (label-#2) at (2.7,1.5) {\bf \Huge #2};%
\end{scope}%
}%
\setlength{\xlen}{5.7cm}%
\setlength{\ylen}{8cm}%
\newcommand{\mzncars}[2]{%
    \foreach \label [count=\n from 0] in {#1} {%
        \ifdefstring{\label}{-}{%
        \phantom{\mzncar[(\n * \xlen, #2 * \ylen)]{\label}}%
        }{%
        \mzncar[(\n * \xlen - 0.3cm, #2 * \ylen)]{\label}%
        }%
        \coordinate (coo-\label) at (\n * \xlen, #2*\ylen);%
        \draw[USred, line width=1mm] (coo-\label) -- +(0.8*\xlen, 0);%
    }%
}%

\tikzset{
    bp-circ/.style={
        circle,
        draw=USred,
        fill=USred,
        inner sep=3pt,
        label={#1},
        solid
    },
    bp-line/.style={
        line width=2.5pt, rounded corners=25pt, USred
    },
  every label/.style={font=\fontsize{28}{48}\selectfont}
}

\begin{tikzpicture}
\normalfont
\def\ysep{3}
\node (right) at (0, 0) {
\begin{tikzpicture}[scale=.5, xscale=-0.85, yscale=0.9]
\mzncars{9, 8, 7, 6, 5, 4, 3, 2, 1}{0}

\coordinate (mid-4) at ($(label-4)+(0, 1.5*\ysep)$);
\coordinate[very thick] (end-4) at ($(label-1) + (0, 1.5*\ysep)$);
\draw[bp-line]  (label-4) -- (mid-4) -- (end-4);
\node[bp-circ, label={left:4}] at (end-4) {};
\coordinate (mid-6) at ($(label-6)+(0, 2.5*\ysep)$);
\coordinate[very thick] (end-6) at ($(label-5) + (0, 2.5*\ysep)$);
\draw[bp-line]  (label-6) -- (mid-6) -- (end-6);
\node[bp-circ, label={left:6}] at (end-6) {};
\coordinate (mid-7) at ($(label-7)+(0, 3.5*\ysep)$);
\coordinate[very thick] (end-7) at ($(label-6) + (0, 3.5*\ysep)$);
\draw[bp-line]  (label-7) -- (mid-7) -- (end-7);
\node[bp-circ, label={left:7}] at (end-7) {};

\coordinate (mid-8) at ($(label-8)+(0, 4.5*\ysep)$);
\coordinate[very thick] (end-8) at ($(label-2) + (0, 4.5*\ysep)$);
\draw[bp-line]  (label-8) -- (mid-8) -- (end-8);
\node[bp-circ, label={left:8}] at (end-8) {};

\coordinate (mid-9) at ($(label-9)+(0, 5.5*\ysep)$);
\coordinate[very thick] (end-9) at ($(label-6) + (0, 5.5*\ysep)$);
\draw[bp-line]  (label-9) -- (mid-9) -- (end-9);
\node[bp-circ, label={left:9}] at (end-9) {};
\end{tikzpicture}
};
\end{tikzpicture}
        }
        \caption{}
    \end{subfigure}
    \caption{The parking blueprint (on the left) and the priority blueprint (on the right) of our running example.}
\end{figure}
\end{ex}

A final observation. the action of the bird`s eye permutation on the priority tree described in Lemma \ref{le:action_weary_permutation} respects the bonsai decomposition. We summarize this in the following lemma.

\begin{lem}
Let $b=(B_1, B_2, \ldots, B_k)$ be the bonsai sequence of $\PT(\pi)$, and let $b_1 = 0$ and  $b_i=|B_1|+|B_2|+\ldots+|B_{i-1}|$ for $i > 1$, then the subtree of $\CT(\pi)$ corresponding to bonsai  $B_i$ is an increasing tree, rooted at $b_{i}+1$, and its nodes  are labelled
$\{b_{i} + 1, b_i + 2, \ldots, b_{i+1}\}.$
\end{lem}
\begin{proof}
This follows from Lemma \ref{le:nodes_parking_tree} recalling that bonsai park one at a time in increasing order of roots.
\end{proof}

Priority trees can also be used to define a new partition for the set of parking functions. 
The equivalence classes of the \emph{priority partition} are  indexed by increasing trees (priority trees), which are well known to be in bijection with permutations. They can also be labeled by subexceedant functions and by priority blueprints.

The priority partition is a refinement of the well-studied partition of the set of parking functions  defined by setting to parking functions that sort into the same sequence to be equivalent. Note that the number of equivalence classes of this second partition is given by the Catalan numbers.

\section{An equidistribution result between statistics for rooted trees and parking functions}
\label{se:equidistribution}

We begin by showing that the weary bijection is record-preserving. Next, we define further interesting statistics and show that they are preserved by the weary map.

\subsection{The weary bijection is record preserving.}

\label{sec:record_bijection} 

We define the set of \emph{records of a parking function}  $\pi$, denoted by \emph{$\Rec(\pi)$}, as the set of records of its bird's eye permutation, so that $\Rec(\pi) = \Rec(\omega_\pi)$. Recall that the records of a permutation are also called its left-to-right maxima.

\begin{thm} 
\label{thm:bijections} Let $\mathcal{T}_n$ be the set of rooted trees of order $n$, and $\mathcal{PF}_n$ be the set of parking functions of length $n$. Then, the  weary bijection is record-preserving.

\end{thm}
 It suffices to show that for any parking function $\pi$,
\[
\Rec(\pi) = \Rec(\PT(\pi)). 
\]

\begin{lem} 
\label{lem:ancestors to the left}
All ancestors in $\PT(\pi)$ of a non-root node  $v \in [n]$   lie to the left of $v$ in the word $\omega_\pi$. Hence, the records  of $\omega_\pi$ are  records of $\PT(\pi)$. 
\end{lem}

\begin{proof}
Given $\omega_i \in [n]$, since the priority vector is subexceedent $\pi(\omega_i) - 1 < i $. On the other hand, the map defining $\PT(\pi)$ verifies 
\[ 
p_{\pi}(\omega_i) = \omega_\pi(\pi(\omega_i) - 1).
\]
This shows that $p_{\pi}(v)$ (this is, the parent of $v$ in $\PT(\pi)$) lies to the left of $v$ in the word $ \omega_\pi$, for any node $v \in [n]$. 
\end{proof}

This result also follows from the recursive description of the process of parking cars on a tree. The record of the $k$-th bonsai, after replacing the label of the root with $\circ$, can be read in $\omega_\pi$ by analyzing the subword found from the $k$-th record up to (but not including) the next one.
As it turns out, the converse of Lemma \ref{lem:ancestors to the left} also holds.
\begin{lem}
\label{le:records_pf}
The records of $\PT(\pi)$ are records of $\omega_\pi$.
\end{lem}
\begin{proof}
    Write $\omega_\pi = \omega_1 \omega_2\cdots \omega_n$ and set $\omega_0 = \circ$.
    Let $\omega_i$ be a record of $\PT(\pi)$.
    The path from $\omega_i$ to the root takes the form 
    \[
    \omega_i, \omega_{j_1}, \dots, \omega_{j_h}=\circ
    \]
    and we know that $i > j_1 > j_2 >\dots > j_h = 0$ by Lemma \ref{lem:ancestors to the left}. 
    
    We argue by contradiction. Suppose that $\omega_i$ is not a left-to-right maxima. Let $\omega_j$ be the right-most element of $\omega_\pi$ to the left of $\omega_i$ that is greater than $\omega_i$. Then, there exists a unique index $k$ such that $j_k 
    >j > j_{k + 1} = \pi(\omega_{j_k}) - 1$. 
    But since $\omega_j$ is parked between $\omega_{j_k}$ and its preferred space $\pi(\omega_{j_k})$, we must have $\omega_j < \omega_{j_k}$. This contradicts our choice of $\omega_j$, since $j < j_k < i$.
\end{proof}


We conclude this work by showing that our record-preserving bijection gives the equidistribution of the following statistics between rooted trees and parking functions. 

\subsection{Statistics of rooted trees}
\label{se:stats_trees}
Let $T$ be a rooted tree of order $n$. We consider the following statistics on $T.$ 
\begin{description}

\item[\emph{Record }]  
 The set of records of a rooted tree $T$ is denoted by \emph{$\Rec(T)$}.

    \item[\emph{Degree of the root }] We denote by \emph{$\deg_\circ(T)$} the number of children of the root in $T$.
    This statistic is the natural analogue for the number of connected components of a rooted forest as removing the root from a rooted tree yields a rooted forest.

    \item[\emph{Children sequence }] 
    The sequence $\memph{\chseq(T)} = (\tau_0, \tau_1 \cdots\tau_{n})$, where $\tau_i$ is the number of nodes of $T$ with exactly $i$ children.
        The sequence $\chseq(T)$ is a composition of $n + 1$ with $n + 1$ parts, some of them being equal to zero. In particular, the number of \emph{leaves} of $T$ is equal to $\tau_0.$

    \item[\emph{Waiting time }]
    Define \emph{$\wait(T)$} as the sum of the total number of cars that park before each node, measured from the moment the node becomes unblocked until the moment it parks, and including the node itself. This quantity can be expressed in terms of the priority tree as $\wait(T) = \diff(\pt(T))$, where $
    \diff(T) = \sum_{i=1}^n (i -p_T(i))
    $. 
 
    \item [\emph{Priority small ascents }] 
    A node $x$ is said to be a \emph{small ascent} if its parent is $x - 1$.
    Define \emph{$\psa(T)$} as the number of small ascents of the priority tree $\pt(T)$. This is, the number of nodes $x$ in $T$ such that the parent $y$ of $x$ is visited just before $x$ in the priority search. 
\end{description}

\subsection{Statistics on parking functions}
\label{Se:stad_parking}
In addition to records, we consider the following statistics for parking functions.  Let $\pi$ be a parking function of length $n.$ 

\begin{description}

\item[\emph{Record }] A record of a parking function $\pi$ is a left-to-right maxima of its bird's eye permutation. The set of record of $\pi$ is denoted by \emph{$\Rec(\pi)$}.

\item[\emph{Ones }] The number of appearances of 1 in $\pi$, denoted by \emph{$\ones(\pi)$}. Note that the only parking function that can have no ones
is the empty parking function. 

\item[\emph{Multiplicity sequence }] The sequence $\memph{\mult(\pi)} = ( \mu_0, \mu_1, \cdots, \mu_{n})$, where $\mu_i$ is the number of elements in $[n + 1]$ that appear with multiplicity $i$ in $\pi$. 
In particular, an element $i$ in $[n + 1]$ is said to be an \emph{absent element} of $\pi$ if the number of appearances of $i$ in $\pi$ is 0. ($n+1$ is always absent). 

\item[\emph{Probes }] The number of attempts (succesful or not) before all cars find their parking spaces, denoted by \emph{$\probes(\pi)$}. 

\item[\emph{Lucky cars }] A car is said to be a \emph{lucky car} of $\pi$ if it parks at its preferred space. The number of lucky cars of $\pi$ is denoted by \emph{$\lucky(\pi)$}.

\end{description}

Note that displacement statistic of \cite{Stanley_parking}, defined as the total number of failed attempts before all cars find their parking spaces, is just a shift of the probes statistic.

\subsection{The equidistribution theorem}
\label{se:thm_equidistribution}
The weary bijection provides a bijective proof of the following equidistribution result.

\begin{thm} \label{thm:statistics}
The sequences of statistics
 \begin{align}  
 \label{eqn:sextuples}
 (\Rec, \wait, \psa, \deg_\circ, \chseq) 
  \ \text{ and } \ 
  (\Rec, \probes, \lucky, \ones, \mult)
\end{align}
are equidistributed when the first sequence is considered over  the second over the set of rooted trees of order $n$, and the second over the set of parking functions of length $n$. 
\end{thm}

Before proceeding to the proof of Theorem \ref{eqn:sextuples}, let us consider some specifications of the equidistributed statistics that are of particular interest. Since the children sequence of the tree $T$ is equal to the multiplicity sequence of the parking function $\pi$, the number of leaves is equal to the number of absent elements of $\pi$ (i.e, the number of elements in $[n + 1]$ that do not appear as the preference of any car).

Moreover, since the number of probes of $\pi$ is equal to the waiting time of $T$, it immediately follows that the lucky cars of $\pi$ (the cars that park in their preferred spot) are exactly those nodes whose label is one larger than their parents' label. That is, the priority small ascents of $T.$ We keep this pair of statistics appearing explictly in statement of Theorem  \ref{thm:statistics} because of their importance in the study of parking functions.

\begin{proof}  There is a single tree in $\mathcal T_0$, namely $\circ$. For this tree, we have 
$
\wait(\circ) = \psa(\circ) = \deg_\circ(\circ) = 0.
$
Similarly, there is a single parking function in $\mathcal{PF}_0$, the empty parking function $\epsilon$, and it verifies
$
\probes(\epsilon) = \lucky(\epsilon) = \ones(\epsilon) = 0.
$
On the other hand, we have $\chseq(\circ) = \mult(\epsilon) = (1, 0, \dots, 0)$ and $\Rec(\circ) = \Rec(\epsilon) = \emptyset$.

Fix $n \geq 1$ and let $T \in \mathcal T_n$. We have already established that $\Rec(T) = \Rec(\pi_T)$.
The number of failed attempts of car $\omega_T(i)$ for the parking function $\pi_T$ is $i - \pv(T)(i) = i - p_{\pt(T)}(i) - 1$. Consequently, the displacement of $\pi_T$ (written $\dis(\pi_T)$) is given by
    \[
    \dis(\pi_T) = \sum_{i = 1}^n (i - p_{\pt(T)}(i)) - n = \wait(T) - n.
    \]
    This shows $\wait(T) = \dis(\pi_T) + n = \probes(\pi_T)$.

    Recall that a car $i$ is lucky for $\pi_T$ if it makes no failed attempts, and by the above $p_{\pt(T)}(i) = i - 1$. Therefore, a car $\omega_\pi(i)$ is lucky if and only if $i$ is a small ascent of $\pt(T)$. This shows $\psa(T)= \lucky(\pi_T)$. 
    
    By definition, $\pi_T(i) = 1$ if and only if $p_T(i)=\circ$. Hence $\deg_\circ(T) = \ones(\pi_T)$. On the other hand, equation \eqref{Def:weary_preference} implies that the number of children of a node $\omega_T(x)$ in $T$ is the number of appearances of $x + 1$ in $\pi_T$.
    From this we deduce that the number of nodes with $k \geq 0$ children in $T$ coincides with the number of symbols in $[n+1]$ that appear in $\pi_T$ with multiplicity $k$, and so the equality $\chseq(T) = \mult(\pi)$ follows.  
\end{proof}

When translated to the language of generating functions, Theorem \ref{thm:statistics} becomes Equation (\ref{eqn: statistics g.f.}):
    \begin{multline}
\label{eqn: statistics g.f.}
  \sum_{T \in \mathcal T} x^{\Rec(T)}  y^{\wait(T)}z^{\psa(T)} t^{\deg_\circ(T)} w^{\chseq(T)} q^{\ord(T)}\\
  = 
  \sum_{\pi \in \mathcal{PF}} x^{\Rec(\pi)} y^{\probes(\pi)} z^{\lucky(\pi)} t^{\ones(\pi)} w^{\mult(\pi)}q^{\len(\pi)}.
    \end{multline}
where $x^{\{a_1, a_2, \dots,a_k\}} = x_{a_1} x_{a_2} \cdots x_{a_k}$ and $w^{(b_1, b_2, \dots, b_k)} = w_{b_1} w_{b_2} \cdots w_{b_k}$,
and where $\ord(T)$ denotes the order of  $T$, and $\len(\pi)$ the length of $\pi$.

\begin{remark}
\label{re:kreweras_enumerator}
    Let $T$ be a rooted tree labelled with $[n]$.
A \emph{tree inversion} is a pair $(i,j)$ such that $i>j$ and
the node $i$ lies on the unique path from the root of $T$ to the node $j$.

    Kreweras \cite{kreweras1980famille} showed that the inversion enumerator on rooted trees coincides with the total displacement enumerator on parking functions. This is, \[
    \sum_{\pi} z^{\dis(\pi)} = \sum_T z^{\inv(T)}.
    \]
    This result has been generalized and extended in subsequent works. For instance,  Stanley and Yin \cite{Stanley_parking} established the equidistribution of a 4-tuple of statistics on trees and parking functions that includes this pair of statistics. 
     Although some statistics may be equidistributed when considered separately, this does not imply that they are jointly equidistributed. Therefore, it is natural to ask whether our statistics can be incorporated to this tuple. 
    
    The 2-tuple of statistics $(\rec, \dis)$ over parking functions and the 2-tuple of statistics $(\rec, \inv)$ over rooted trees
    are \emph{not} equidistributed. Indeed, for $n = 2$ we have 
    \begin{align*}
   & \sum_{\pi} z^{\dis(\pi)}t^{\rec(\pi)} = zt^2 + t^2 + t,
   & \text{ whereas }
   & \sum_T z^{\inv(T)}t^{\rec(T)} = 2t^2 + zt. 
    \end{align*}
    The explicit data behind this claim is presented in Tables \ref{tab:pfs} and \ref{tab:trees}.

Furthermore, none of the pairs of statistics studied by Stanley and Yin are jointly equidistributed with the record statistics. Consequently, both equidistribution results are independent;  see (\ref{final:StanleyYin}) of the final comments. 

\begin{table}[ht]
\centering
\normalfont
\renewcommand{\arraystretch}{1.2}

\newlength{\tableW}
\setlength{\tableW}{0.42\textwidth}

\begin{minipage}{0.48\textwidth}
\centering
\begin{tabularx}{\tableW}{c|*3{>{\centering\arraybackslash}X}}
\toprule
$\pi$        & 11 & 12 & 21 \\
$\omega_\pi$ & 12 & 12 & 21 \\
\midrule
$\rec(\pi)$  & 2  & 2  & 1  \\ 
$\dis(\pi)$  & 1  & 0  & 0  \\ 
\bottomrule
\end{tabularx}
\captionof{table}{Records and displacement in parking functions of length $2$.}
\label{tab:pfs}
\end{minipage}
\begin{minipage}{0.48\textwidth}
\centering
\begin{tabularx}{\tableW}{c|*3{>{\centering\arraybackslash}X}}
\toprule
\multirow{2}{*}{$T$} 
& \multirow{2}{*}{\hspace*{0.2cm}\raisebox{-.3\height}{\begin{tikzpicture}[
    level distance=0.3cm, 
    level 1/.style={sibling distance=0.5cm},
    level 3/.style={sibling distance=1cm},
    level 4/.style={sibling distance=0.6cm},
    tree/.style={
        anchor=north, 
    },
    record/.style={
        circle,
        draw=black,
        fill=black,
        inner sep=1.2pt,
        label={#1},
        solid
    },
    record/.default={},
    non-record/.style={
        circle,
        draw = black,
        fill = white,
        inner sep=1.2pt,
        solid
    },
]
\node[non-record, label={right:$\circ$}] (N-0) {}
child {	node[record, label={right:1}] (N-1) {}
	child {	node[record, label={right:2}] (N-2) {}}};    
\end{tikzpicture}}} 
& \multirow{2}{*}{\hspace*{0.2cm}\raisebox{-.3\height}{\begin{tikzpicture}[
    level distance=0.3cm, 
    level 1/.style={sibling distance=0.5cm},
    level 3/.style={sibling distance=1cm},
    level 4/.style={sibling distance=0.6cm},
    tree/.style={
        anchor=north, 
    },
    record/.style={
        circle,
        draw=black,
        fill=black,
        inner sep=1.2pt,
        label={#1},
        solid
    },
    record/.default={},
    non-record/.style={
        circle,
        draw = black,
        fill = white,
        inner sep=1.2pt,
        solid
    },
]

\node[non-record, label={right:$\circ$}] (N-0) {}
child {	node[record, label={right:2}] (N-2) {}
	child {	node[non-record, label={right:1}] (N-1) {}}};   
\end{tikzpicture}}} 
& \multirow{2}{*}{\hspace*{0.2cm}\raisebox{-.3\height}{\begin{tikzpicture}[
    level distance=0.35cm, 
    level 1/.style={sibling distance=0.5cm},
    level 3/.style={sibling distance=1cm},
    level 4/.style={sibling distance=0.6cm},
    tree/.style={
        anchor=north, 
    },
    record/.style={
        circle,
        draw=black,
        fill=black,
        inner sep=1.2pt,
        label={#1},
        solid
    },
    record/.default={},
    non-record/.style={
        circle,
        draw = black,
        fill = white,
        inner sep=1.2pt,
        solid
    },
]

\node[non-record, label={right:$\circ$}] (N-0) {}
child {	node[record, label={right:1}] (N-1) {}}
child {	node[record, label={right:2}] (N-2) {}};
    
\end{tikzpicture}}} \\
\vphantom{$\omega_\pi$} & & & \\ 
\midrule
$\rec(T)$ & 2 & 1 & 2 \\ 
$\inv(T)$ & 0 & 1 & 0 \\ 
\bottomrule
\end{tabularx}
\captionof{table}{Records and inversions in rooted trees of order $2$.}
\label{tab:trees}
\end{minipage}
\end{table}

\end{remark}

\subsection{Record duality.}

Let $A$ be a set of rooted trees, $B$  a set of parking functions.
We say that a $A$ and $B$ are in \emph{record duality} if the image of $A$  under the weary bijection is $B$.

The duality between increasing trees and subexceedant functions is the central theme of this work.
 More examples of record duality are easy to come by:  Path graphs with permutations.
Catalan trees (rooted trees in which any node is either a leaf or has two children) and parking functions such that any symbol appears with multiplicity 0 or 2, 
$k$-ary trees (rooted trees in which every node has at most $k$ children) and parking functions such that any symbol appears with multiplicity at most $k$.  Indeed, any family of rooted trees, defined exclusively by a set of conditions on its children sequence, will be in duality with the family of parking functions, with that same set of conditions in its multiplicity sequence. 
Alternatively, other statistics may also be considered. For example, we could  fix a  set of records, or a bound on the number of lucky cars or ones or probes, or a long etc.
In all these cases, the corresponding statistics are equidistributed.

\begin{cor} \label{thm:statistics_dual} 
Let $A$ be a set of rooted trees, and $B$  a set of parking functions. If $A$ and $B$ are in record duality, then the  statistics 
 \begin{align*}  
 (\Rec, \wait, \psa, \deg_\circ, \chseq) \ \text{ and } \ 
& (\Rec, \probes, \lucky, \ones, \mult).
\end{align*}
are equidistributed when the first sequence of statistics is considered over the set of  rooted trees in $A$ and the second  over the set of parking functions  in $B$.
\end{cor}

Let $A$ be a set of rooted trees, and $B$  a set of parking functions. At the level of generating functions, Corollary \ref{thm:statistics_dual} says that if $A$ and $B$ are in record duality, then
     \begin{multline*}
  \sum_{T \in A }  
  x^{\Rec(T)}
  y^{\wait(T)}z^{\psa(T)} t^{\deg_\circ(T)} w^{\chseq(T)} q^{\ord(T)}
   \\
  = 
  \sum_{\pi \in B} x^{\Rec(\pi)}y^{\probes(\pi)} z^{\lucky(\pi)} t^{\ones(\pi)} w^{\mult(\pi)}q^{\len(\pi)}.
    \end{multline*}


\section{Two enumerative questions}
\label{se:hook}

We consider the following two enumerative questions. On the one hand, we fix a permutation $\omega$ in $\SS_n$ and ask: How many rooted trees/parking functions have $\omega$ as its bird eye permutation? On the other hand, we fix a increasing tree $T$  as ask: How many rooted trees/parking functions have $T$ as its priority tree?  

\subsection{The number of parking functions with a fixed bird's eye permutation.}
Fix  $\omega$ in $\SS_n$, and  $k$ in $[n]$. The \emph{priority run} of $k$ is the longest run of $\omega$ ending at $k$ with all its elements $\le k$. We denote the length of the priority run of $k$ by \emph{$\lpr(k)$}. 

\begin{lem} Fix a permutation $\omega$ in $\SS_n$.
Then, the number of rooted trees with priority transversal $\omega$, and the number of parking functions with bird's eye permutation $\omega$ are both equal to
\begin{align}
\label{KW_formula}
\prod_{k=1}^n \lpr(k).
\end{align}
\end{lem}
 
\begin{proof}   Lemma $\ref{lem: omega_T = omega_pi_T}$ implies that the two numbers coincide. Thus, it suffices to prove one of these claims. 
We work in  the setting of parking functions.
Fix $k$ in $[n]$, and let  $\ell= \lpr(k)$. Let $p_1, p_2, \ldots, p_\ell$, with $p_\ell=k$ be the priority run of $k$. By definition of parking functions, this implies that spots $p_1, p_2, \ldots, p_\ell$ are already filled at the moment when car $k$ attempts to park. So, if the preference of $k$ if $p_1$, $p_2$, $\ldots$ or   $p_\ell$, car  $k$ will end up parked at  position, $\omega(k)$. 

On the other hand, the car that parks at the spot to the left of the one containing $p_1$ arrives after $k$, (after all, it has a bigger label). Thus,  when car $k$ attempts to park this spot is empty. If car $k$ would had preferred to park either there or at an earlier spot, it would always find a free spot before  $\omega(i)$.
\end{proof}

\begin{ex} Let $\omega = \pfunction{5, 2, 4, 7, 1, 3, 6, 8, 9}$ be the bird's eye permutation of our running examples. \\

\begin{center}
\normalfont
\renewcommand{\arraystretch}{1.15}
\begin{tabular}{c|*{9}{>{\centering\arraybackslash}p{1.8em}}@{\hspace{0.8em}}}
\toprule
$\omega(i)$
  & 5 & 2 & 4 & 7 & 1 & 3 & 6 & 8 & 9 \\
\midrule
priority run of $i$
  & 5 & 2 & {\small 24} & 7 & 1 & {\small 13} & {\small 136} & {\small 1368} & {\small 13689} \\
$\lpr(i)$
  & 1 & 1 & 2 & 3 & 1 & 2 & 3 & 4 & 5 \\
\bottomrule
\end{tabular}
\end{center}

\smallskip

Therefore, there are 720 parking functions with $\omega$ as its bird eye's permutation, and as many trees with priority transversal $\omega$.
\end{ex}
This result can be traced back to the pioneering work of Konheim and Weiss \cite{KW}. In their paper, Eq.~(\ref{KW_formula}) was used  to show recursively that the number of parking functions of length 
$n$ is is $(n+1)^{n-1}$. It does not seem to exist in the literature in the context of trees.


\subsection{On the number of parking functions with a fixed priority tree.} A hook-length formula of Gessel and Seo (Lemma 2.1 of \cite{Gessel-Seo}) allows us to answer this second question.
Let $T$ be a rooted tree, and let $\pi_T$ be its parking function. We introduce a partial order on $[n],$ called the \emph{parking order} of $T$  on $[n]$, where we say that 
$a \lessdot_T b$ (this is, $b$ covers $a$) if the arc of $b$ in the priority blueprint of $\pi$ lies immediately above the parking spot $a$.  Observe that the weary bijection allows us to also refer to $<_T$ as the parking order of $\pi_T$. 

\begin{figure}[h!]
    \centering
    \normalfont

    \begin{subfigure}{0.18\textwidth}
        \centering
        \resizebox{0.45\textwidth}{!}{\begin{tikzpicture}[
    level distance=0.6cm, 
    level 1/.style={sibling distance=0.7cm}, 
    level 2/.style={sibling distance=0.6cm},
    tree/.style={
        anchor=north, 
    }
]

\node[non-record, label={right:$\circ$}] (T1) {}
    child {	
	    node[record, label={right:1}] (1-4-1) {}
        child{
            node[record, label={right:2}] {}
        }
        }
    child {	
	    node[record, label={right:3}] (1-4-3) {}
	    child {	
		    node[record, label={right:4}] (1-4-2) {}}};
\end{tikzpicture}}
        \caption{}
        \label{subfig:small_priority_tree}
    \end{subfigure}
    \qquad
    \begin{subfigure}{0.28\textwidth}
        \centering
        \resizebox{0.7\textwidth}{!}{\newcommand{\mzncar}[2][(0,0)]{%
\begin{scope}[shift={#1}]%
\draw[draw=black, fill=black, rounded corners=1.2ex, very thick] (0,.5) -- ++(0,1) -- ++(1,0.3) -- ++(3,0) -- ++(1.1,0.05) -- ++(0,-1.3) -- cycle;%
\draw[very thick, rounded corners=0.5ex, fill=black, thick] (1,1.74) -- ++(1,0.7) -- ++(1.6,0) -- ++(0.9,-0.7) -- cycle;%
\draw[fill=white, draw=white, thick] (1.25,.5) circle[radius=.55];%
\draw[fill=white, draw=white, thick] (4,.5) circle[radius=.55];%
\draw[fill=black,semithick] (1.25,.5) circle[radius=.35];%
\draw[fill=black,semithick] (4,.5) circle[radius=.35];%
\node[white] (label-#2) at (2.7,1.5) {\bf \Huge #2};%
\end{scope}%
}%
\setlength{\xlen}{5.7cm}%
\setlength{\ylen}{8cm}%
\newcommand{\mzncars}[2]{%
    \foreach \label [count=\n from 0] in {#1} {%
        \ifdefstring{\label}{-}{%
        \phantom{\mzncar[(\n * \xlen, #2 * \ylen)]{\label}}%
        }{%
        \mzncar[(\n * \xlen - 0.3cm, #2 * \ylen)]{\label}%
        }%
        \coordinate (coo-\label) at (\n * \xlen, #2*\ylen);%
        \draw[USred, line width=1mm] (coo-\label) -- +(0.8*\xlen, 0);%
    }%
}%

\tikzset{
    bp-circ/.style={
        circle,
        draw=USred,
        fill=USred,
        inner sep=3pt,
        label={#1},
        solid
    },
    bp-line/.style={
        line width=2.5pt, rounded corners=25pt, USred
    },
  every label/.style={font=\fontsize{28}{48}\selectfont}
}

\begin{tikzpicture}
\normalfont
\def\ysep{3}
\node (right) at (0, 0) {
\begin{tikzpicture}[scale=.5, xscale=-0.85, yscale=0.9]
\mzncars{4, 3, 2, 1}{0}

\coordinate (mid-3) at ($(label-3)+(0, 1.5*\ysep)$);
\coordinate[very thick] (end-3) at ($(label-1) + (0, 1.5*\ysep)$);
\draw[bp-line]  (label-3) -- (mid-3) -- (end-3);
\node[bp-circ, label={left:3}] at (end-3) {};
\end{tikzpicture}
};
\end{tikzpicture}}
        \caption{}
        \label{subfig:small_blueprint}
    \end{subfigure}
    \qquad
    \begin{subfigure}{0.22\textwidth}
        \centering
        \resizebox{0.55\textwidth}{!}{\begin{tikzpicture}[
    level distance=0.6cm, 
    level 1/.style={sibling distance=0.7cm}, 
    level 2/.style={sibling distance=0.6cm},
    tree/.style={
        anchor=north, 
    }
]

\node[non-record, label={right:$\circ$}] (T1) {}
    child {	
	    node[record, label={right:3}] (1-4-1) {}
        child{
            node[non-record, label={right:1}] {}
        }
        child{
            node[non-record, label={right:2}] {}
        }
        }
    child {	
	    node[record, label={right:4}] (1-4-3) {}};
\end{tikzpicture}}
        \caption{}
        \label{subfig:small_cover_tree}
    \end{subfigure}

    \vspace{0.5cm}

    \begin{subfigure}{\textwidth}
        \centering
        \resizebox{\textwidth}{!}{%
            \usetikzlibrary{matrix}
\begin{tikzpicture}[
    level distance=0.5cm, 
    level 1/.style={sibling distance=0.7cm}, 
    level 2/.style={sibling distance=0.6cm},
    tree/.style={
        anchor=north, 
    }
]
\def\codeheight{15pt}

\newcommand{\basetree}{
    
}

\matrix[matrix of nodes, row sep=0.5cm, column sep=1cm, nodes={tree}, ampersand replacement=\&] {
    \node[non-record, label={right:$\circ$}] (T1) {}
    child {	
	    node[record, label={right:1}] (1-4-1) {}
        child{
            node[record, label={right:2}] {}
        }
        }
    child {	
	    node[record, label={right:3}] (1-4-3) {}
	    child {	
		    node[record, label={right:4}] (1-4-2) {}}};\&
    \node[non-record, label={right:$\circ$}] (T2) {}
    child {	
	    node[record, label={right:2}] (1-4-1) {}
        child{
            node[non-record, label={right:1}] {}
        }
        }
    child {	
	    node[record, label={right:3}] (1-4-3) {}
	    child {	
		    node[record, label={right:4}] (1-4-2) {}}};\&
    \node[non-record, label={right:$\circ$}] (T3) {}
    child {	
	    node[record, label={right:3}] (1-4-1) {}
        child{
            node[non-record, label={right:1}] {}
        }
        }
    child {	
	    node[record, label={right:4}] (1-4-3) {}
	    child {	
		    node[non-record, label={right:2}] (1-4-2) {}}};\&
    \node[non-record, label={right:$\circ$}] (T4) {}
    child {	
	    node[record, label={right:2}] (1-4-1) {}
        child{
            node[non-record, label={right:1}] {}
        }
        }
    child {	
	    node[record, label={right:4}] (1-4-3) {}
	    child {	
		    node[non-record, label={right:3}] (1-4-2) {}}};\&
    \node[non-record, label={right:$\circ$}] (T5) {}
    child {	
	    node[record, label={right:1}] (1-4-1) {}
        child{
            node[record, label={right:3}] {}
        }
        }
    child {	
	    node[record, label={right:4}] (1-4-3) {}
	    child {	
		    node[non-record, label={right:2}] (1-4-2) {}}};\&
    \node[non-record, label={right:$\circ$}] (T6) {}
    child {	
	    node[record, label={right:2}] (1-4-1) {}
        child{
            node[record, label={right:3}] {}
        }
        }
    child {	
	    node[record, label={right:4}] (1-4-3) {}
	    child {	
		    node[non-record, label={right:1}] (1-4-2) {}}};\&
    \node[non-record, label={right:$\circ$}] (T7) {}
    child {	
	    node[record, label={right:1}] (1-4-1) {}
        child{
            node[record, label={right:2}] {}
        }
        }
    child {	
	    node[record, label={right:4}] (1-4-3) {}
	    child {	
		    node[non-record, label={right:3}] (1-4-2) {}}};\&
    \node[non-record, label={right:$\circ$}] (T8) {}
    child {	
	    node[record, label={right:3}] (1-4-1) {}
        child{
            node[record, label={right:2}] {}
        }
        }
    child {	
	    node[record, label={right:4}] (1-4-3) {}
	    child {	
		    node[non-record, label={right:1}] (1-4-2) {}}};\&
            \\};
\def\pfsep{2}
\node (pf1) at ($(T1) + (0, -\pfsep)$) {\footnotesize $\pfunction{1,2,1,4}$}; 
\node (pf2) at ($(T2) + (0, -\pfsep)$) {\footnotesize $\pfunction{2,1,1,4}$}; 
\node (pf3) at ($(T3) + (0, -\pfsep)$) {\footnotesize $\pfunction{2,4,1,1}$}; 
\node (pf4) at ($(T4) + (0, -\pfsep)$) {\footnotesize $\pfunction{2,1,4,1}$}; 
\node (pf5) at ($(T5) + (0, -\pfsep)$) {\footnotesize $\pfunction{1,4,2,1}$}; 
\node (pf6) at ($(T6) + (0, -\pfsep)$) {\footnotesize $\pfunction{4,1,2,1}$}; 
\node (pf7) at ($(T7) + (0, -\pfsep)$) {\footnotesize $\pfunction{1,2,4,1}$}; 
\node (pf8) at ($(T8) + (0, -\pfsep)$) {\footnotesize $\pfunction{4,2,1,1}$};

\end{tikzpicture}
        }
        \caption{}
    \label{subfig:trees_with_same_priority_tree}
    \end{subfigure}

    \caption{(a) The parking tree $P$ of $\pi = \pfunction{1,2,1,4}$. (b) The parking blueprint of $\pi$. (c)  The hasse diagram of the parking order of $\pi$. (d) The 8 trees with priority tree $P$ and the 8 parking functions with priority vector $\pi$.}
    \label{fig:same priority tree}
\end{figure}
\begin{lem}[Lemma 2.1, \cite{Gessel-Seo})]
\label{le:Gessel-Seo}
Fix an increasing tree $T$ of order $n$. Let $H$ be the Hasse diagram of the parking order $<_T$.
Then, the number of rooted trees with priority tree $T$ and the number of parking functions with priority tree $T$  are both equal to  
\begin{align}
\label{Eq:GesseL_Seo}
          \frac{n!}{\prod_{v\in H} h(v)}
\end{align}
\end{lem}

\begin{proof} Let  $\pi_T$ be the parking function of $T$. 
Let $a, b$ and $c$ be nodes of the Hasse diagram (cars trying to park).
First observe that   if $a<_Tb$ in the parking order, then $a<b$. 
On the other hand, if $a$  and $b$ both cover $c$ and $a<b$, then $a<_Tb$. These two observations follow from the definition of parking blueprint. 
Therefore, if we root each connected component of the Hasse diagram  at it maximum element, we end up with a rooted forest, that we may turn into a rooted tree by adding an extra node $\circ$ which we attach all roots of it. The result is always a decreasing tree.

At this point the result then follows immediately from the Gessel-Seo formula (Lemma 2.1 of \cite{Gessel-Seo}).
As a final comment, we can ignore the extra node $\circ$ in the Gessel-Seo formula as taking it into account is equivalent to multiplying and dividing by $n$ (because $h(\circ)=n$).
\end{proof}

\begin{ex} Consider the subexceedant function $\pfunction{1, 2, 1, 4}$. Construct its priority tree $P$ (Figure \ref{subfig:small_priority_tree}), its parking blueprint (Figure \ref{subfig:small_blueprint}), and the Hasse diagram of its parking order (Figure \ref{subfig:small_cover_tree}). The hook-length formula tells us that there are $8$ trees with priority tree $P$. 
On the other hand, let $T$ be the increasing tree of our running example. (Example \ref{ex_parkingtree}). The hook-length formula tells us that there are 210 trees with priority tree $T$.
\end{ex}

Observe that
Lemma \ref{le:Gessel-Seo} determines the cardinalities of equivalence classes of the priority partition.

\begin{figure}[t]
    \centering
    \normalfont
    
    \begin{subfigure}{0.42\textwidth}
        
        \raisebox{5mm}{
            \resizebox{0.9\textwidth}{!}{%
                \newcommand{\mzncar}[2][(0,0)]{%
\begin{scope}[shift={#1}]%
\draw[draw=black, fill=black, rounded corners=1.2ex, very thick] (0,.5) -- ++(0,1) -- ++(1,0.3) -- ++(3,0) -- ++(1.1,0.05) -- ++(0,-1.3) -- cycle;%
\draw[very thick, rounded corners=0.5ex, fill=black, thick] (1,1.74) -- ++(1,0.7) -- ++(1.6,0) -- ++(0.9,-0.7) -- cycle;%
\draw[fill=white, draw=white, thick] (1.25,.5) circle[radius=.55];%
\draw[fill=white, draw=white, thick] (4,.5) circle[radius=.55];%
\draw[fill=black,semithick] (1.25,.5) circle[radius=.35];%
\draw[fill=black,semithick] (4,.5) circle[radius=.35];%
\node[white] (label-#2) at (2.7,1.5) {\bf \Huge #2};%
\end{scope}%
}%
\setlength{\xlen}{5.7cm}%
\setlength{\ylen}{8cm}%
\newcommand{\mzncars}[2]{%
    \foreach \label [count=\n from 0] in {#1} {%
        \ifdefstring{\label}{-}{%
        \phantom{\mzncar[(\n * \xlen, #2 * \ylen)]{\label}}%
        }{%
        \mzncar[(\n * \xlen - 0.3cm, #2 * \ylen)]{\label}%
        }%
        \coordinate (coo-\label) at (\n * \xlen, #2*\ylen);%
        \draw[USred, line width=1mm] (coo-\label) -- +(0.8*\xlen, 0);%
    }%
}%

\tikzset{
    bp-circ/.style={
        circle,
        draw=USred,
        fill=USred,
        inner sep=3pt,
        label={#1},
        solid
    },
    bp-line/.style={
        line width=2.5pt, rounded corners=25pt, USred
    },
  every label/.style={font=\fontsize{28}{48}\selectfont}
}

\begin{tikzpicture}
\normalfont
\def\ysep{3}
\node (right) at (0, 0) {
\begin{tikzpicture}[scale=.5, xscale=-0.85, yscale=0.9]
\mzncars{9, 8, 7, 6, 5, 4, 3, 2, 1}{0}

\coordinate (mid-4) at ($(label-4)+(0, 1.5*\ysep)$);
\coordinate[very thick] (end-4) at ($(label-1) + (0, 1.5*\ysep)$);
\draw[bp-line]  (label-4) -- (mid-4) -- (end-4);
\node[bp-circ, label={left:4}] at (end-4) {};
\coordinate (mid-6) at ($(label-6)+(0, 2.5*\ysep)$);
\coordinate[very thick] (end-6) at ($(label-5) + (0, 2.5*\ysep)$);
\draw[bp-line]  (label-6) -- (mid-6) -- (end-6);
\node[bp-circ, label={left:6}] at (end-6) {};
\coordinate (mid-7) at ($(label-7)+(0, 3.5*\ysep)$);
\coordinate[very thick] (end-7) at ($(label-6) + (0, 3.5*\ysep)$);
\draw[bp-line]  (label-7) -- (mid-7) -- (end-7);
\node[bp-circ, label={left:7}] at (end-7) {};

\coordinate (mid-8) at ($(label-8)+(0, 4.5*\ysep)$);
\coordinate[very thick] (end-8) at ($(label-2) + (0, 4.5*\ysep)$);
\draw[bp-line]  (label-8) -- (mid-8) -- (end-8);
\node[bp-circ, label={left:8}] at (end-8) {};

\coordinate (mid-9) at ($(label-9)+(0, 5.5*\ysep)$);
\coordinate[very thick] (end-9) at ($(label-6) + (0, 5.5*\ysep)$);
\draw[bp-line]  (label-9) -- (mid-9) -- (end-9);
\node[bp-circ, label={left:9}] at (end-9) {};
\end{tikzpicture}
};
\end{tikzpicture}
            }}%

        \caption{}
    \end{subfigure}
    \qquad
    \begin{subfigure}{0.34\textwidth}
        \centering
        \raisebox{0mm}{
            \begin{tikzpicture}[
    level distance=0.6cm, 
    level 2/.style={sibling distance=1cm},
    level 3/.style={sibling distance=1cm},
    level 4/.style={sibling distance=0.6cm},
    tree/.style={
        anchor=north, 
    }
]

\node[non-record, label={right:$\circ$}] (N-circ) {}

child{node[record, label={right:9}] (N-9) {}
child {	node[non-record, label={right:8}] (N-8) {}
	child {	node[non-record, label={left:7}] (N-6) {}
		child {	node[non-record, label={left:6}] (N-3) {}
			child {	node[non-record, label={left:5}] (N-1) {}}}}
	child {	node[non-record, label={right:4}] (N-7) {}
		child {	node[non-record, label={right:3}] (N-2) {}}
		child {	node[non-record, label={right:2}] (N-4) {}}
		child {	node[non-record, label={right:1}] (N-5) {}}}}};
    
\end{tikzpicture}%
        }
        
        \caption{}
    \end{subfigure}
    \caption{(a) The priority blueprint of the tree $T$ of Example \ref{ex:running_as_image_of_tree}.
    (b) The Hasse diagram of the parking order of $T$. }
\end{figure}


\section{Final Comments} \label{sec:final-comments}
\begin{enumerate}[wide, labelindent=8pt]

\item 
All combinatorial structures introduced in this work, as well as the main bijection, have been implemented in Sage  \cite{sage}. The source code can be found in \cite{WearySage}.

\newcommand{\tikzxmark}{%
\tikz[scale=0.23] {
    \draw[line width=0.7,line cap=round] (0,0) to [bend left=6] (1,1);
    \draw[line width=0.7,line cap=round] (0.2,0.95) to [bend right=3] (0.8,0.05);
}}
\newcommand{\tikzcmark}{%
\tikz[scale=0.23] {
    \draw[line width=0.7,line cap=round] (0.25,0) to [bend left=10] (1,1);
    \draw[line width=0.8,line cap=round] (0,0.35) to [bend right=1] (0.23,0);
}}

\item 
\label{final:StanleyYin}
Richard Stanley and  Mei Yin  \cite{Stanley_parking} established an equidistribution result between rooted trees and parking functions somewhat orthogonal to ours.  
However, none of their statistics are jointly equidistributed with the record statistic $\Rec$, not even with the coarser statistic $\rec$. 
 Indeed, observe that the following statistics are not equidistributed when restricted to parking functions and rooted trees with exactly $n$ records (namely, subexceedant functions and increasing trees): 

\begin{itemize}[wide, labelindent=12pt]

    \item[\tikzxmark]  $(\text{nld}, \text{unl})$. The number of non-leaders of any increasing tree is $0$, yet a wide range of values for the unluckyness statistic can be found among the set of subexceedant functions. From $0$, as in $12\cdots n$, up to $n - 1$, as in $1111\cdots1$. 
    \item[\tikzxmark]  $(\text{inv}, \text{dis})$. The number of inversions of an increasing tree is 0. Nontheless, the displacement statistic of subexceedant functions ranges from 0 (as in $123\cdots n$) up to $\binom{n}{2}$ (as in $111\cdots 1$). 
    \item[\tikzxmark]  $(\leaves, \text{des})$. Since the bird's eye permutation of any subexceedant function is the identity, it has 0 descents. The number of leaves in an increasing tree, however, ranges from 1 (in the increasing path graph) up to $n - 1$ (in the star tree).
    \item[\tikzxmark] $(\deg_\circ, \text{rlm})$. The number of right-left maxima of the bird's eye permutation of a subexceedant function is $1$. However, the root of an increasing tree can have any number of children. 
\end{itemize}

\item As stated in the introduction (Section \ref{sec:intro}), the priority-first search we have described is closely related to Prim's algorithm, which takes as input an edge-weighted graph and outputs a spanning tree of minimal weight (where the weight of the spanning tree is defined to be the sum of its edge weights). There, one starts with an arbitrary node of the graph as the only visited node. At each step one follows the smallest weighted edge $uv$ for which the node $u$ is visited and $v$ is unvisited in order to visit the next unvisited node $v$. This process matches the priority traversal if one simply ``pushes up'' the vertex labels of each non-root node and views them as edge weights.

\item An analysis of the  statistics introduced in this paper, in the spirit of work of Diaconis and Hicks \cite{DiaconisHicks} will be very welcomed. 
 It would also be interesting to characterize explicitly further examples pairs of sets of rooted trees/parking functions in record-duality.

\end{enumerate}

\section*{Acknowledgements}
We benefited greatly from discussions with Matthieu Josuat-Vergès, who introduced us to the Gessel–Seo hook formula. We also thank Persi Diaconis, Dominique Foata, and Ira Gessel for their encouragement and support. Finally, we thank the anonymous referees for their valuable comments.

    Authors AL and MR have been partially supported by Grants PID2020-117843GB-
I00 and PID2024-157173NB-I00 funded by MCIN/AEI/10.13039/501100011033 and by FEDER, UE. Author ST is supported by \href{10.3030/101070558}{FoQaCiA} which is funded by the European Union and NSERC.
 The second author would like to thank the IRIF of the Université de Paris Cité for their hospitality, and  the opportunity to conduct part of this research at their facilities. 

\bibliographystyle{halpha-abbrv} 
\bibliography{references} 
\end{document}